\documentclass[trsc,nonblindrev]{informs3-neutral} 
\OneAndAHalfSpacedXI 

\usepackage{natbib}
 \bibpunct[, ]{(}{)}{,}{a}{}{,}%
 \def\bibfont{\small}%
\TheoremsNumberedThrough     

\EquationsNumberedThrough    

\usepackage{afterpage}
\usepackage{color}
\usepackage{enumitem}
\usepackage{setspace}
\usepackage{algorithmic}
\usepackage{algorithm}
\usepackage[table]{xcolor}
\usepackage{multirow}
\usepackage{microtype}
\usepackage[hidelinks]{hyperref}
\usepackage{booktabs}
\usepackage{rotating}

\usepackage{commands}

\newcommand{\myred}[1]{#1}

\newcommand{\myblue}[1]{#1}
\renewcommand{\myblue}[1]{#1} 

\newenvironment{add}
{}
{}

\newenvironment{del}
{}
{}

\title{Electric Vehicle Fleets: Scalable Route and Recharge Scheduling through Column Generation}

\RUNAUTHOR{Axel Parmentier, Rafael Martinelli, Thibaut Vidal}
\RUNTITLE{Electric Vehicle Fleets: Scalable Route and Recharge Scheduling}

\MANUSCRIPTNO{TS-2021-0139}

\begin{document}

\begin{center}

\begin{LARGE}
Electric Vehicle Fleets: Scalable Route and  \vspace*{0.3cm} \\ Recharge Scheduling through Column Generation
\end{LARGE}

\vspace*{01cm}

\textbf{Axel Parmentier$^1$, Rafael Martinelli$^2$, Thibaut Vidal$^{3,4}$} \\
$^1$ \textsc{Cermics}, \'Ecole des Ponts, Marne la Vall\'ee, France \\
$^2$ Department of Industrial Engineering, Pontifical Catholic University of Rio de Janeiro, Brazil\\
$^3$ CIRRELT \& SCALE-AI Chair in Data-Driven Supply Chains, Department of Mathematics and Industrial Engineering, Polytechnique Montréal, Canada\\
$^4$ Department of Computer Science,  Pontifical Catholic University of Rio de Janeiro, Brazil\\
\url{axel.parmentier@enpc.fr} \hfill \url{martinelli@puc-rio.br} \hfill
\url{thibaut.vidal@polymtl.ca}
\vspace*{0.8cm}

\end{center}

\noindent
\textbf{Abstract.} 
The rise of battery-powered vehicles has led to many new technical and methodological hurdles. Among these, the efficient planning of an electric fleet to fulfill passenger transportation requests still represents a major challenge.
This is due to the specific constraints of electric vehicles, bound by their \myblue{battery} autonomy and necessity of recharge planning, and the large scale of \myblue{the operations}, which challenges existing optimization algorithms. The purpose of this paper is to introduce a scalable column generation approach for \myblue{routing and scheduling in this context. Our algorithm relies on four main ingredients: 
(i) a multigraph reformulation of the problem based on a characterization of non-dominated charging arcs,
(ii) an efficient bi-directional pricing algorithm using tight backward bounds,}
(iii) sparsification approaches permitting to decrease the size of the subjacent graphs dramatically, and 
(iv) a diving heuristic, which locates near-optimal solutions in a fraction of the time needed for a complete branch-and-price. Through extensive computational experiments, we demonstrate that our approach significantly outperforms previous algorithms for this setting, leading to accurate solutions for problems counting several hundreds of requests.
\vspace*{0.2cm}

\noindent
\textbf{Keywords.} Routing and scheduling, Electric vehicles, Column generation, Diving heuristics.

\vspace*{0.5cm}

\thispagestyle{empty}
\pagenumbering{arabic}

\section{Introduction}
\label{sec:introduction}

To cut down pollution and carbon footprint, major cities worldwide have established regulations to progressively phase out conventional internal combustion engine (ICE) vehicles and replace them with subsidized electric vehicles (EV). Public bus transportation has already shifted towards battery-powered vehicles in many cities \citep{Li2016c,Mahmoud2016,Wang2017c}.
In contrast, taxi companies and mobility-on-demand services take more time to transition as they often depend on the driver's willingness to switch to EVs. Still, EV taxi fleets are progressively becoming more widespread in various cities \citep[see, e.g.,][]{Dunne2017,DeJong2018,London2018,Hu2018,Scorrano2020}. \myblue{Finally, freight transportation operators increasingly integrate electric vehicles into their fleets \citep{BloombergNEF2021}. This ongoing transformation has led to new market segments focused on emission-free transportation, and various incentives have been established to provide advantageous leasing conditions, dedicated priority queues, and charging infrastructures for EV drivers.}

\myblue{Despite all the actions taken, electric fleets remain subject to different and challenging operational constraints. Whereas buses operate on a fixed itinerary and require fast-charging (or battery-swap) stations in a limited number of locations, passenger transportation operators servicing a varying set of addresses require scattered charging stations to operate efficiently. Competition is also intense in the transportation sector, such that efficient planning of customer-to-vehicle assignments, dead-heading trips, and recharging actions is critical for customer satisfaction and profitability.}
When the charging infrastructure is established, the optimization of vehicle itineraries and recharging trips to meet the demand for a set of timetabled trips can be formalized as an electric vehicle scheduling problem (EVSP --  \citealt{Bodin1983,Wen2016}). It is common to register hundreds of transportation or service requests in densely populated areas, leading to planning problems of large size that call for new solution~paradigms.

This paper contributes to overcoming the existing methodological gaps and the need for efficient and scalable solution algorithms. 
\myblue{It first formalizes the concept of non-dominated charging arcs, considering linear and non-linear charging, leading to an efficient algorithm for enumerating these arcs and to a compact definition of the EVSP on a multigraph.} Then, it introduces exact and heuristic column generation algorithms grounded on the monoid pricing paradigm by \cite{Parmentier2018}. Our pricing algorithm \myblue{exploits high-quality backward bounds on path resources to reduce the number of paths. To further enhance efficiency on large instances, we develop sparsification techniques to reduce the size of the multigraph.} Then, we embed the column generation in a specialized diving algorithm that uses strong diving up to a given depth. As seen in our experiments, these techniques permit us to efficiently solve large-scale problems with several hundreds of requests. The contributions of this paper therefore~are:
\begin{enumerate}
    \item \myblue{A characterization of non-dominated charging arcs, considering possible sequences of charging stations along with linear or non-linear charging rates.}
	\item A new pricing algorithm for EVSPs, leading to the first mathematical programming algorithm able to solve all existing 100-customers instances for the problem.
	\item Tailored sparsification and diving strategies, which permits scaling up the algorithm to cases counting half a thousand visits and finding near-optimal solutions in those cases.
	\item An extensive numerical campaign, which evaluates the impact of critical methodological components and demonstrates the method's applicability on a wide range of instances. \myred{All the numerical experiments can be reproduced using our open-source code located at \url{https://github.com/axelparmentier/ElectricalVSP-ColumnGeneration}}.
	\item More generally, this work paves the way towards a new generation of scalable and flexible algorithms for EVSPs.
\end{enumerate}

The remainder of this paper is as follows. Section~\ref{sec:problem_statement_and_related_works} formally defines the EVSP and briefly reviews related works. \myblue{Next, Section~\ref{sec:StationSequenceScheduling} explains how non-dominated charging arcs corresponding to station sequences can be identified and characterized by a charging function in the general non-linear case as well as in the particular case of linear charging.
The proposed pricing algorithm is described in Section \ref{sec:solution_methods}, along with the exact branch-and-price approach and driving heuristic built upon it.} Section~\ref{sec:numerical_results} describes our computational experiments, and Section~\ref{sec:conclusion} finally concludes.

\section{Problem Statement and related works}
\label{sec:problem_statement_and_related_works}

\vspace*{1em}
\myblue{
\noindent
\textbf{Setting.}
In a similar fashion as in \cite{Wen2016}, the EVSP can be defined as follows. Let $[0,\tm]$ be a planning horizon, where $0$ and $\tm \in \bbR$ respectively correspond to the beginning and the end of the day.}
Let~$\Tr$ be a set of timetabled passenger trips, also called \emph{services}, that a transportation \myblue{(e.g., taxi, bus, or shuttle)} company should operate during the day.
Let $\De$ be the set of depots where the vehicles start and end, and let~$\St$ be a set of stations where vehicles can charge their battery.
Each service $\tr \in \Tr$ is planned to start at~$\tb[\tr] \in \bbR$ and end at $\te[\tr] \in \bbR$.
\myblue{We use the convention~$\tb[\tr] = \tm$ and $\te[\tr] = 0$ for each depot $\tr$ in $\calD$.
We further assume that} the number of vehicles stationed at each depot is not limited, but that each vehicle is characterized by a fixed use cost $\cv$ and equipped with a battery of capacity~$\mc$.
For notation convenience, we introduce the set $\Ac = \Tr \cup \De \cup \St$ and set $\fu[\ac] = 0$ for $\ac \in \De \cup \St$.
For each pair $\ac[1],\ac[2]$ of elements of~$\Ac$,
we denote respectively by $\ti[\ac_1,\ac_2]$, and $\co[\ac[1],\ac[2]]$ the driving time, and cost required to go from the location of (the end of) $\ac[1]$ to the location of (the beginning~of)~$\ac[2]$.

\vspace*{1em}
\myblue{
\noindent\textbf{Batteries linear discharge and non-linear recharge.}
Battery discharge is typically linear, i.e., energy consumption while driving typically does not depend on the current charge level.
Each service $\tr$ consumes a quantity $\fu[\tr]$ of energy, and we denote by $\fu[\ac[1],\ac[2]]$ the energy required to go from the location of $\ac[1]$ to the location of $\ac[2]$.
}

\myblue{
Batteries can be recharged at any charging station. 
Typically, the charge level follows an ordinary differential equation $\frac{d\ell}{dt} = i^\mathrm{ch}(\ell)$, where $i^\mathrm{ch}$ is a non-increasing non-negative function. Let us denote $\chf : [0,\mc]\times \bbR_+ \rightarrow [0,\mc]$ the flow of this ordinary differential equation.
If a vehicle with an initial battery level $\ell$ spends a duration $\tau$ in a station, then it leaves the station with battery level $\chf(\ell,\tau)$. Since $\varphi$ is the flow of an ordinary differential equation, and given the properties of $i^\mathrm{ch}$, $\ell,\tau \mapsto \chf(\ell,\tau)$ is non-decreasing, $\tau \mapsto \chf(\ell,\tau)$ is concave, $\ell \mapsto \varphi(\ell,\tau) - \ell$ is non-increasing, and $\chf(\ell,\tau+\tau') = \chf(\chf(\ell,\tau),\tau')$.
Given $\ell$ and $\ell'$ in $[0,\mc]$, let us also denote by $\tch(\ell,\ell')$ the charge time needed to go from $\ell$ to $\ell'$
$$\tch(\ell,\ell')=
\begin{cases}
     \min \{\tau\colon \varphi(\ell,\tau)=\ell'\}, & \text{if } \ell \leq \ell', \\
     - \tch(\ell',\ell), & \text{otherwise.}
\end{cases}
$$
\begin{remark}
    \myred{Charge is supposed to be identical in all charging stations. 
    If this hypothesis is not satisfied, the algorithms in this paper still work if vehicles can visit only one station between two tasks. Otherwise, some dominance results exploited in algorithms (e.g., for  Proposition~\ref{prop:optimalSchedulingOfStationSequence}) do not hold.}
\end{remark}

}

\myblue{
\noindent\textbf{Multigraph representation and general objective.}
For each depot $\de$ in $\De$, we introduce the directed multigraph $D_{\de} = (V_{\de},A_{\de})$ defined as follows. 
The vertex set $V = \calT \cup \{o,d\}$ contains the tasks~$\calT$ as well as dummy vertices $o$ and $d$ corresponding to route start and end at the depot~$\de$.
If there is enough time between the end of $u$ and the beginning of $v$ to travel from $u$ to $v$, i.e., if
$\te[v] + \ti[u,v] \leq \tb[v]$
we have a \emph{non-charging} $u$-$v$ arc $a$ in $A_{\de}$ that models that a vehicle can travel directly from the end of $u$ to the beginning of~$v$. If there is enough time between the end of $u$ and the beginning of $v$, the vehicle can stop in one or several stations in $\St$ and charge its battery. Some sequences of stops in charging stations are irrelevant. We introduce in Section~\ref{sec:StationSequenceScheduling} a notion of dominance between station sequences as well as efficient algorithms to preprocess all relevant sequences of stops. At this point, let us simply suppose that we have a set of \emph{charging arcs} $a$ between $u$ and $v$ that model the different relevant sequence of one or several stops in charging stations that the vehicle can do when it travels from $u$ to $v$. We denote by $A_{\de}^{\mathrm{ch}} \subseteq A_{\de}$ the set of charging arcs.
Since there can be an $u$-$v$ arc only if \myred{$v$} begins after \myred{$u$} ends, $D_{\de}$ is acyclic.
In equations, we sometimes use notation $a=(u,v)$ to mean that $u$ is the tail of $a$ and $v$ its head.

We assume that each arc $a$ in $A$ can be characterized by a cost $c_a$ and a mapping $\fc{a} : [0,M] \cup \{-\infty\} \rightarrow [0,M] \cup \{-\infty\}$ that gives the battery level $\fc{a}(\lin)$ at the end of $a$ if the battery level is $\lin$ at the beginning of $a$. \myred{Note that $\fc{a}$ is not a vector. The arrow in the notation will be used to distinguish forward-backward computations of battery level.} Level $-\infty$ is used to identify that the vehicle runs out of battery. If $a$ is the non-charging arc between $u$ and $v$, then we have:
\begin{equation}
    \label{eq:nonChargingArcsCostandLevel}
    c_a = \co[u,v] + \co[v]
    \quad \text{and} \quad
     \fc{a}(\lin) = \pos(\lin - \fu[u,v] - \fu[v])
    \quad \text{where} \quad
    \pos(\ell) =
    \begin{cases}
        \ell & \text{if }\ell \geq 0\\
        - \infty & \text{otherwise.}
    \end{cases}\vspace*{0.4cm}
\end{equation}
Section~\ref{sec:StationSequenceScheduling} provides the characterization of $c_a$ and $\fc{a}$ for charging arcs, which is more involved.}

\myblue{
A vehicle route $r$ is 
an $o$-$d$ path $o,a_1,v_1,\ldots,v_k,a_{k+1},d$  in one of the $D_{\de}$.
We denote by $\tr \in \ro$ and $a\in\ro$ the fact that the service $\tr$ and the station sequence $a$ belong to route $\ro$, respectively.
The \emph{cost} $c_r$ of route $r$ is the sum of a fixed vehicle cost $\cv$ and the cost of its arcs:
$$c_r = \cv + \sum_{a \in \ro}c_a.$$
Route $r$ is \emph{feasible} if it can be operated without running out of battery, that is, if and only if
$\fc{a_{k+1}} \circ \fc{a_k}  \circ \cdots \circ \fc{a_1}(\mc) \geq 0$, \myred{where $\circ$ denotes function composition.}
We denote by $\fR$ the set of feasible routes.
}

With these definitions, the goal of the EVSP is to find a set $\sol \subset \fR$ of feasible routes covering each service and such that the sum of the costs of the routes in $\sol$ is minimum. Given~$x_{\ro}$, a binary variable indicating if a route $\ro$ is in the solution, the problem can be stated as:
\begin{subequations}\label{eq:masterPb}
	\begin{alignat}{2}
	\min_{x} & \sum_{\ro \in \fR} \co[r]x_{\ro}\\
	\mathrm{s.t.} & \sum_{\ro \in \fR\colon \tr \in r} x_{r} = 1 &\quad& \forall \tr \in \Tr \label{eq:partconstraint} \\
	& x_r \in \{0,1\} && \forall r \in \fR.
	\end{alignat}
\end{subequations}

Naturally, the size of $\fR$ grows exponentially with the number of services, making any direct enumeration impractical in most practical cases. We will therefore rely on column generation \citep{Barnhart1998} to circumvent this issue.

\vspace*{1em}
\noindent
\textbf{Related Works.} 
Early studies on vehicle scheduling problems (VSP) date back to \cite{Bodin1983} and often arose from the air transportation literature. The term ``vehicle scheduling'' has been coined to describe routing problems in the presence of determined start dates for the services, therefore establishing a total order. Solution algorithms for VSPs are typically more scalable than their vehicle routing problem (VRP) counterparts since they are not dependent on cycle elimination. \myblue{A similar situation happens when solving vehicle routing problems with time windows when services (e.g., scheduled taxi trips) have a longer duration than the width of their windows of acceptable start time, as this characteristic makes it impossible to have cycles.} Consequently, the canonical VSP with a single depot and no other constraint than the customers' service times is known to be polynomially solvable as a minimum-cost flow problem in an acyclic digraph. However, many immediate extensions of the problem, e.g., with multiple depots or duration constraints on the routes, belong to the NP-hard class.
Among all variants surveyed in \cite{Bodin1983} and \cite{Bunte2009}, the VSP with ``length of path considerations'' (i.e., duration limits) is particularly relevant to our problem. Duration constraints model range-limitations which are typical in electric vehicles, but do not allow for possible recharging. This variant has been solved by \cite{Haghani2002}, \cite{Ribeiro1994}, and \cite{Desrosiers1995} with heuristics and column-generation algorithms.

More recently, VSPs with duration constraints and possible recharging stops have regained attention due to their economic significance in electric mobility systems. \cite{Adler2017} considered the VSP with alternative-fuel vehicles and multiple depots, in which a limited number of nodes act as charging stations. This model only allows full charging and supposes that the charging time is constant, regardless of the energy level upon arrival to the station. A similar model is studied in \cite{Li2014}. Constant-time charging is relevant when considering battery-swapping technologies, but only represents a coarse approximation of mobility systems based on charging technology. 
Therefore, \cite{Wen2016} define the electric VRP (EVSP), which integrates the possibility of partial charging in a time which grows linearly with the amount of electricity recharged. Finally, \cite{VanKootenNiekerk2017} propose a more accurate and general charging model, considering non-linear charging rates as well as possible battery-deterioration effects. Due to the inherent complexity of this charging scheme, discretization techniques are often used for the more complex models.

A similar progression is observable when surveying the literature about electric vehicle routing problems (EVRP). As illustrated in the surveys of \cite{Brandstatter2016,Pelletier2016,Schiffer2018} and \cite{Vidal2020}, after a first wave of studies focused on full constant-time charging, research have progressed towards more sophisticated and realistic problem settings, e.g., with mixed vehicle fleets \citep{Hiermann2016,Hiermann2019}, richer delivery networks \citep{Schiffer2017,Breunig2018}, or considering heterogeneous charging infrastructures \citep{Felipe2014,Keskin2018} and non-linear charging rates \citep{Montoya2017,Pelletier2017}. However, the resulting models are highly complex to solve, such that the vast majority of these articles propose metaheuristic approaches. Moreover, solution methods are tested on benchmark instances that rarely exceed a hundred customers.

Specific to exact methods, with the exception of  \cite{Wen2016}, which is based on a mixed integer programming (MIP) formulation and branch-and-cut algorithm, all state-of-art algorithms for EVSP variants \citep{Li2014,Adler2017,VanKootenNiekerk2017} rely on a set partitioning (or set covering) formulation, and therefore repeatedly solve a pricing problem which can be formulated as a variant of resource-constrained shortest path problem (RCSP) \citep{Irnich2005}. When only full constant-time charging is allowed, the pricing problem is a particular case of the weight-constrained shortest path with replenishment arcs, studied in \cite{Smith2012} and \cite{Bolivar2014}. More realistic charging assumptions (e.g., partial and non-linear charging) lead to more complex RCSPs. In most state-of-the-art EVSP and EVRP column generation algorithms \citep{Desaulniers2016}, pricing problems are solved with labeling techniques, using resource extension functions. Bi-directional search strategies have also been exploited in different algorithms \citep{Tilk2017,Thomas2019,Cabrera2020}. However, despite all the progress done, no existing algorithm for the EVSP can consistently solve instances with more than a hundred customers. To fill this gap, we explore new bi-directional pricing and bounding strategies grounded on the work of \citet{Parmentier2018}. In addition, to scale up the solution approach to large-scale problems, we design efficient sparsification techniques and diving heuristics, allowing a fast heuristic column generation with only a minor loss of solution quality.
\myblue{Our methodology will be described in the next sections, starting with the characterization of the charging arcs, and following with the solution algorithms.}

\section{\myblue{Charging arcs}}
\label{sec:StationSequenceScheduling}

The purpose of this section is to give a formal definition of charging arcs $a$ along with a characterization of their cost $c_a$ and mapping $\fc{a}$.
We also characterize the mapping $\bc{a} \colon \lout \rightarrow \bc{a}(\lout)$ which associates to a target charge level $\lout$ the minimum charge $\bc{a}(\lout)$ needed at the beginning of $a$ to reach the end of $a$ with charge level at least $\lout$. This mapping $\bc{a}$ is leveraged by our solution algorithms.
All these results rely on the notion of station sequence, which we now introduce.

\subsection{\myblue{Station sequence}}
\myblue{
When traveling from the depot to a service, between two services, or from a service to the depot, a vehicle can possibly stop in one or more charging stations.
Let $u$ and $v$ in $\De \cup \Tr$ be such that it is possible to operate $v$ after $u$, that is, $\te[u] + \ti[u,v] \leq \tb[v]$.
A \emph{station sequence $a$ between $u$ and $v$} in $\calV$ is a sequence $u = s_0,s_1,\ldots,s_k,s_{k+1}=v$ with $s_1,\ldots,s_k$ in $\calS$ such that:
\begin{equation}\label{eq:stationSequenceTime}
    \te[u] + \sum_{i=1}^{k+1} \delta_{s_{i-1},s_i} \leq \tb[v], \quad \text{and} \quad
    \begin{cases}
    e_{u} + e_{u,s_1} \leq \mc, \\
    e_{s_{i-1},s_i} \leq \mc  & \text{ for $i$ in $\{1,\ldots,k\}$} \\
    e_{s_k,v} + e_v \leq \mc. 
    \end{cases}
\end{equation}
These constraints mean that there must be enough time between the end of $v$ and the beginning of $u$ to travel along the station sequence, and that, given two successive elements of the sequence, we can reach the second if we have a full battery at the beginning of the first.
The cost of a station sequence $a$ is:
\begin{equation}\label{eq:stationSequenceCost}
    \co[a] = \sum_{i=1}^{k+1}\Big(\co[s_{i-1},s_i]\Big) + \co[v].
\end{equation}
Contrarily to services, the arrival time $\tb[s]$ and the departure time $\te[s]$ from a station $s$ are not fixed. 
A scheduling $\btr[a]$ of $a$ is a vector $(\tb[s_1],\te[s_1],\ldots, \tb[s_k],\te[s_k])$ giving the departure time and arrival time for each station $s$ in $a$ such that:}
\footnote{\myblue{With this convention, no time is spent waiting in the sequence. This is done w.l.o.g. since recharging is superior to waiting at a station in all aspects in the current problem statement.}}
\begin{equation}
\label{eq:stationSequenceTimes}
\myblue{
 \tb[s_i] = \te[s_{i-1}] + \ti[s_{i-1},s_i] \text{ for } i \in \{1,\ldots,k+1\}  \quad \text{and} \quad 
 \tb[s_i] \leq \te[s_i] \text{ for } i \in \{1,\ldots,k\}.}
\end{equation}
\myblue{If we have a scheduling $\btr[a]$ and the charge level $\lin$ at the end of $u$, we can deduce the charge level at any point during the station sequence. 
For $i \in \{1,\ldots,k+1\}$,
let $\levb[s_i][a](\lin,\btr[a])$ be the charge level at beginning of~$s_i$. 
And for $i \in \{0,\ldots,k+1\}$, let $\levb[s_i][a](\lin,\btr[a])$ be the charge level at the end of~$s_i$. We have the following recursive definition:
\begin{equation}
\label{eq:batteryLevel}
\begin{array}{l}
    \levb[s_i][a](\ib,\btr[a]) = 
    \begin{cases}
        \ib, & \text{if } i=0, \\
        - \infty, & \text{if } 1 \leq i \leq k \text{ and } \te[s_i] < \tb[s_i], \\
        \varphi(\levb[s_i][a](\ib,\btr[a]), \te[s_i] - \tb[s_i]),  & \text{if } 1 \leq i \leq k \text{ and } \te[s_i] \geq \tb[s_i], \\
        \pos(\levb[s_i][a](\ib,\btr[a]) - e_{s_i}), & \text{if } i = k+1,
    \end{cases}\vspace*{0.4cm}
    \\
     \text{and} \quad
    \levb[s_i][a](\lin,\btr[a]) =  \pos(\lev[s_{i-1}][a](\ib,\btr[a]) - \fu[s_{i-1},s_i]).
    \quad \text{Recall that }
    \pos(\ell) =
    \begin{cases}
        \ell & \text{if }\ell \geq 0,\\
        - \infty & \text{otherwise.}
    \end{cases}\vspace*{0.4cm}
\end{array}
\end{equation}
With this definition, $\lev[s_i][a](\ib,\btr[a]) = -\infty$ if the battery level becomes negative at some point.
}
\myblue{
Since $\btr[a]$ does not impact the costs, the vehicle should always use the scheduling $\btr[a]$ that leads to the maximum charge level at the end of the sequence. 
Due to the non-linear recharge, the relevance of a scheduling $\btr[a]$ may depend on the initial charge level $\lin$. 
Let $\fc{a}(\lin)$ be the charge at the end of $a$ under an optimal scheduling given the initial charge level $\lin$ at the end of $u$, i.e.,
\begin{equation}
    \label{eq:optimalScheduilngProblem}
    \fc{a}(\lin) = \max \big\{\lev[v][a](\ib,\btr[a])\colon \btr[a] \text{ satisfies \eqref{eq:stationSequenceTimes}}\big\}.
\end{equation}
\noindent
\textbf
{Station sequences and arcs.}
We use the same notation $a$ for arcs of the multigraph $D_{\de}$ and stations sequences because arcs are going to be defined as special cases of station sequences.
Non-charging and charging $u$-$v$ arcs are station sequences $u = s_0,s_1,\ldots,s_k,s_{k+1}=v$ with $k=0$ and $k>0$, respectively.
Remark that the characterization of $c_a$ and $\fc{a}$ for non-charging arcs in Equation~\eqref{eq:nonChargingArcsCostandLevel} corresponds to the definitions in Equations~\eqref{eq:stationSequenceCost} and~\eqref{eq:optimalScheduilngProblem} with $k=0$.
}

\subsection{\myblue{Optimal Schedule Characterization}}

\myblue{
Solving~\eqref{eq:optimalScheduilngProblem} may seem difficult a priori. Yet, a simple recursion enables solving it in the general case. Let $a$ be the station sequence $u=s_0,s_1,\ldots,s_{k},s_{k+1}=v$. 
A larger initial charge is always better. Indeed, since $\ell \mapsto \varphi(\ell,\tau)$ is non-decreasing, $\lin  \mapsto \lev[s_i][a](\ib,\btr[a])$ is non-decreasing for any $s_i$ in $a$ and $\btr[a]$, which gives that $\lin \mapsto \fc{a}(\lin)$ is non-decreasing. Next, because recharging with $\varphi$ is more efficient when the charge level is low, we will formally show that spending the smallest feasible amount of time in all stations but the last is optimal. This leads us to the following theorem.}

\myblue{
\begin{theorem}
\label{prop:optimalSchedulingOfStationSequence}
Consider the schedule $\btr[a,*]$ defined recursively as follows:\vspace*{-0.3cm}
\begin{align*}
    &\tb[s_i][a,*](\ib)= \te[s_{i-1}][a,*](\lin) + \ti[s_{i-1},s_i] \\
    &\te[s_i][a,*](\lin)=
	\begin{cases}
	\tb[s_i][a,*](\lin) +
	\left\{
	\begin{array}{ll}
	   0 & \text{if }  \lev[s_i][a]\big(\btr[a,*](\lin)\big) \geq \fu[s_i,\ac[i+1]]  \\
	   \tch\big(\lev[s_i][a]\big(\btr[a,*](\lin)\big), \fu[s_i,s_{i+1}]\big)& \text{otherwise.}
	\end{array}
	\right. 
	& \text{if } 1 \leq i \leq k-1 \\
	\tb[v] - \ti[s_k,v]
	& \text{if } i=k.
	\end{cases}
\end{align*}
If $\btr[a,*](\lin)$ is feasible, i.e., $\tb[s_k][a,*](\lin) \leq \te[s_k][a,*](\lin)$, then it is an optimal solution of Problem~\eqref{eq:optimalScheduilngProblem}. Otherwise, $\fc{a}(\lin) = - \infty$.
\end{theorem}}

\myblue{
\noindent Let us denote $\levb[s_i][a,*](\lin)$ and $\lev[s_i][a,*](\lin)$ by $\levb[s_i][a](\lin,\btr[a,*])$ and $\lev[s_i][a](\lin,\btr[a,*])$, respectively.
}

\vspace*{1em}
\myblue{
\noindent \textbf{Discussion.}
In practice, an optimal schedule $\btr[a]$ as given in Theorem~\ref{prop:optimalSchedulingOfStationSequence} can be risky as it typically suggests arriving in a station with an empty battery. For practical applications, we recommend a leveled version where a station is reached with a given reserve battery level $M_0$. From a mathematical viewpoint, this is equivalent to replacing $M$ with $M-M_0$ in the problem definition. Note that introducing setup times or costs, if needed, for recharging does not change the structure of optimal schedules, as these characteristics can be directly integrated into the definition of the travel times and travel costs to reach the stations.
}

\subsection{\myblue{Station sequences dominance and formal definition of $A_{\de}$}}

\myblue{
Consider two station sequences $a$ and $a'$ between $u$ and $v$. We say that 
$$\text{$a$ \emph{dominates}~$a'$ \quad if} \quad \co[a] \leq \co[a'] \quad \text{and} \quad \fc{a} \geq \fc{a'}.$$
If $a$ dominates $a'$, then given any feasible route $r'$ that contains $a'$, the route $r$ that is obtained by replacing $a$ by $a'$ in $r$ is also feasible, by monotonicity of the $\fc{a_i}$ for $a_i$ in $r'$, and with a non-greater cost $c_{r} \leq c_{r'}$. Therefore, we can limit our search to a subset of optimal solutions that dominates all the other station sequences. 
We can now formally define the set of arcs $A_{\de}$ introduced in Section~\ref{sec:problem_statement_and_related_works} as a \emph{set of non-dominated station sequences}, i.e.,
a set $A$ of station sequences such that, (1) for any $u$, $v$, and station sequence $a$ between $u$ and $v$, either $a$ is in $A$ or there is a stations sequence $a'$ between $u$ and $v$ in $A$ such that $a'$ dominates $a$, and (2) no $u$-$v$ arc in $A$ is dominated by another $u$-$v$ arc in $A$.

\vspace*{1em}
\noindent \textbf{Practical computation of $A_{\de}$.}
Let us now describe how to build a set of non-dominated station sequences between $u$ and $v$. 
Let $S = \big([0,M] \cup \{-\infty\}\big) \times \bbR$,
and $\preceq$ be the homogeneous binary relation on $S \times S$ defined by 
$$ 
\begin{cases}
    (-\infty,t) \preceq (\ell',t'), & \text{for any }(\ell',t') \in S \quad \text{ and $t$ in $\bbR$}, \\ 
    (\ell,t) \preceq (-\infty,t')  \quad \text{if and only if} \quad \ell = -\infty, & \text{for any $t,t'$ in $\bbR$}, \\
    (\ell,t) \preceq (\ell',t') \quad \text{if and only if} \quad 
    t + \tch(\ell,\ell') \geq t', & \text{for $\ell,\ell' \in [0,\mc]$ and $t,t'$ in $\bbR$}.
\end{cases}
$$
The relation $\preceq$ permits comparing the charge level at different time points. 
Given $t\leq t'$ and $(\ell,t) \preceq (\ell',t')$, it can be interpreted as follows: If we have a battery level of $\ell$ at $t$, and we spend all the time between $t$ and $t'$ charging the vehicle, then the charge level at $t'$ is less than or equal to $\ell'$.
We prove in Appendix~B that $\preceq$ is a preorder: it is reflexive because $\varphi(\ell,0) = \ell$, and transitive because $\varphi(\ell,\tau + \tau') = \varphi\big(\varphi(\ell, \tau),\tau'\big)$ for $\tau,\tau'\geq0$.

\begin{theorem}\label{prop:stationSequenceDominance}
    Let $a$ and $a'$ be the two station sequences $u=s_0,s_1,\ldots,s_k, s_{k+1} = v$ and $u=s'_0,s'_1,\ldots,s'_{k'}, s'_{k'+1} = v$ between $u$ and $v$.
    \begin{enumerate}
        \item If $a$ and $a'$ coincide up to one station, i.e., if there exists $i$ such that $s_j = s'_j$ for all $j\leq i$, then $\tb[s_i][a,*](\lin) = \tb[s_i][a',*](\lin)$ and $\levb[s_i][a,*](\lin) = \levb[s_i][a',*](\lin)$ for any $\lin$. 
        \item If $a$ and $a'$ coincide starting from one station, i.e., if there exists $i$ and $i'$ such that $k-i = k'-i'$ and $s_{i + j} = s'_{i' + j}$ for $j$ in $\{0,\ldots,k-i\}$, and if:
        $$\sum_{j=1}^i \co[s_{j-1},s_j] \leq \sum_{j=1}^{i'} \co[s'_{j-1},s'_j] 
        \quad \text{and} \quad
        \big(\levb[s_i][a,*],\tb[s_i][a,*](\lin)\big) 
        \succeq 
        \big(\levb[s'_{i'}][a',*],\tb[s'_{i'}][a,*](\lin)\big) \text{ for all } \lin \in [0,\mc],$$
        then $a$ dominates $a'$.
    \end{enumerate}
\end{theorem}
Theorem~\ref{prop:stationSequenceDominance} shows that the problem of building $A_{\de}$ amounts to generating the Pareto front of a multiobjective shortest path problem. 
Indeed, consider the digraph $D_{u,v}$ whose vertex set is $\calS \cup \{u,v\}$ and whose arc set is:
$$\big\{(u,s)\colon s\in \St\big\} \cup \big\{(s,s')\colon s,s'\in \St\big\} \cup  \big\{(s,v) \colon s \in \St\big\}. $$
Given a $u$-$s$ path $P$ of the form $u,s_0,s_1,\ldots,s_i=s$, we define $c_
P$ as $\sum_{j = 1}^i \co[s_{i-1},s_i]$.
The first point of Theorem~\ref{prop:stationSequenceDominance} ensures that we can define the arrival time $\tb[][P](\lin)$ and the charge level $\levb[][P](\lin)$ of a $u$-$s$ path $P$ in $D_{u,v}$  as the arrival time $\tb[s][a,*](\lin)$ in $s$ and the charge level $\levb[s][a,*](\lin)$ at the arrival in $s$ for any station sequence $a$ that starts by $P$.
Given two $u$-$s$ path $P$ and $P'$ in $D_{u,v}$, we say that $P$ dominates $P'$ if:
$$c_P \leq c_{P'} \quad \text{and} \quad \big(\levb[][P](\lin), \tb[][P](\lin) \big) \succeq \big(\levb[][P'](\lin), \tb[][P'](\lin) \big) \text{ for all } \lin \in [0,\mc].$$
The second point of Theorem~\ref{prop:stationSequenceDominance} shows that building the non-dominated station sequences between $u$ and $v$ amounts to building the non-dominated paths starting in $u$ in $D_{u,v}$.
It can be done with any algorithm generating the Pareto front of non-dominated paths (e.g., the one of~\citealt{kergosien2021efficient}).
Since the number of stations is relatively limited,
the resulting problems can be efficiently solved \myred{provided that the mappings $\lin \mapsto \levb[][P'](\lin)$ can be encoded efficiently and dominance checks are efficient}. 
As seen in our numerical experiments, the time needed by this preprocessing is only a tiny fraction of the total computing time.
}

\subsection{\myblue{Forward-backward check of route feasibility}}
\myblue{
Because $\lin \mapsto \fc{a}(\lin)$ is non-decreasing, we can define an inverse as:
\begin{equation}\label{eq:definition_of_bca}
    \bc{a}(\lout) = \min \big\{\lin \colon \fc{a}(\lin) \geq \lout \big\}.
\end{equation}
Its output $\bc{a}(\lout)$ can be interpreted as the minimum charge needed at the beginning of $a$ to be able to reach the end of $a$ with a charge level $\lout$.}
\myblue{These two functions 
$$
    \begin{array}{rcl}
         \fc{a}\colon \underline{L} & \rightarrow &  \underline{L}\\
        \lin & \mapsto &   \fc{a}(\lin)
    \end{array}
    \quad \text{and} \quad
    \begin{array}{rcl}
        \bc{a}\colon \overline{L} & \rightarrow &  \overline{L}\\
        \lout & \mapsto &  \bc{a}(\lout)
    \end{array}
    \quad\text {where}\quad
    \begin{cases}
    \overline{L} &=[0,M] \cup \{+ \infty\} \\
    \underline{L} &= [0,M] \cup \{- \infty\} 
    \end{cases}
$$
enable to check the feasibility of a route $r$ in a forward-backward way.
}
\myblue{
\begin{lemma}
\label{lem:bidirectionalRouteFeasibility}
Let $r=v_0,a_1,\ldots,a_{k+1},v_{k+1}$ be a route and $i \in \{0,\ldots,k+1\}$.
Then $r$ is feasible if and only if 
$\fc{a_i}\circ \ldots \circ \fc{a_1}(\mc) \geq \bc{a_{i+1}}\circ\ldots \circ \bc{a_{k+1}}(0)$.
\end{lemma}

\textsc{Proof.}
Using reverse induction on $i$, we get that 
$\bc{a_{i+1}}\circ\ldots \circ \bc{a_{k+1}}(0) =  \min \big\{\lin \colon \fc{a_{k+1}}\circ \ldots \circ \fc{a_{i+1}}(\lin) \geq 0 \big\}$, which gives us the following result.
Hence $\fc{a_{k+1}}\circ \ldots \circ \fc{a_1}(\mc)$ is non-negative if $\fc{a_i}\circ \ldots \circ \fc{a_1}(\mc) \geq \bc{a_{i+1}}\circ\ldots \circ \bc{a_{k+1}}(0)$, and equal to $-\infty$ otherwise.
\halmos

\vspace*{1em}
Lemma~\ref{lem:bidirectionalRouteFeasibility} states that $r$ is feasible if and only if the charge level at the end of stations sequence $a_i$ is sufficient to operate the rest of the route. This permits to jointly exploit forward and backward recursion in pricing algorithms.}

\vspace*{1em}
\noindent
\textbf{Keys for an efficient implementation.}
When applied to challenging instances, the pricing algorithm calls millions of times the routines computing $\fc{a}$ and $\bc{a}$. Therefore, we need an efficient algorithm to compute them. In the general non-linear case, 
\myred{we expect the station sequences to be generally short (only occasionally two stations or more in practice). Indeed, traveling to a station is costly, so it is rarely interesting to make detours through several stations.} Hence, $\fc{a}(\lin)$ can be efficiently computed using the formula for $\btr[a,*]$. \myred{Alternatively, if time and charge level are discretized, we could tabulate these functions.}
The backward computation of $\bc{a}(\lout)$ is slightly more involved. Let $a$ be the arc corresponding to the station sequence $u=s_0,s_1,\ldots,s_k,s_{k+1} = v$ between $u$ and $v$, and consider $\lout \in [0,\mc]$.
Let $\lin = \bc{a}(\lout)$.
Suppose that $\lin < + \infty$.
Then Theorem~\ref{prop:optimalSchedulingOfStationSequence} ensures the existence of $i^*$ be such that:
\vspace*{1em}
\begin{itemize}
    \item $\tb[s_i][a,*](\lin) = \te[v] + \sum_{j=1}^{i}\ti_{s_{j-1},s_j}$ and $\lev[s_{i-1}][a](\lin,\btr[a,*](\lin) > \fu[s_{i-1},s_i] $ for $i\leq i^*$, 
    \item $\tb[s_i][a,*](\lin) > \te[v] + \sum_{j=1}^{i}\ti_{s_{i-1},s_i}$ and $\lev[s_{i-1}][a](\lin,\btr[a,*](\lin) =  \begin{cases}
    \fu[s_{i-1},s_i] & \text{if } i< k \\
    \lout + e_{s_k,v} + e_{v} & \text{if } i= k 
    \end{cases} $ for $i > i^*$.\vspace*{1em}
\end{itemize}
In particular, if $i> i^*$, then the charge level at the arrival in $s_i$ is equal to $0$. 
This gives a simple backward algorithm to detect $i^*$ by constructing $\btr[a,*](\lin)$ in a backward way. 
If $\lout > \mc - e_{s_{k},v} - e_v$, it means $\bc{a}(\lout) = +\infty$. 
Otherwise, for $i$ from $k$ to $1$ and while $i^*$ has not been identified, do the following operations.  Compute $t' = \tb[s_{i+1}] - \delta_{s_i,s_{i+1}} - \tch(0,l)$, where $l=e_{s_i,s_{i+1}}$ if $i<k$ and $l =\lout + \fu[s_k,v] + \fu[v]$ otherwise. 
If $t' > \te[v] + \sum_{j=1}^{i}\ti_{s_{i-1},s_i}$, then $i^* < i$, and we set $\tb[s_i] = t'$ and we can continue the backward algorithm.
On the contrary, if $t' \leq \te[v] + \sum_{j=1}^{i}\ti_{s_{i-1},s_i}$, it means that $i=i^*$. 
Therefore, we know that $\tb[s_i][*]= \te[v] + \sum_{j=1}^{i}\ti_{s_{i-1},s_i}$ and $\te[s_i][*]= \tb[s_{i+1}] - \delta_{s_i,s_{i+1}}$.
We then solve for $\tilde l$ in the equation
$ \varphi\big(\tilde l, \te[s_i][*] - \tb[s_i][*]\big) = l$,
where $\tilde l$ is going to be the charge level at the arrival in $s_i$.
In the most general case, this can be done by dichotomy.
We can then compute $l_0 = \tilde l + \sum_{j=1}^i \fu[s_{i-1},s_{i}]$.
If $l_0 > \mc$, then $\bc{a}(\lout)= + \infty$. 
Otherwise, $\bc{a}(\lout) = l_0$.
Finally, if we have identified that $i^*=0$, we have $\bc{a}(\lout) = \fu[u,s_1]$.

\subsection{\myblue{Practical Computations in the Linear Case}}

\myblue{
Recharge is linear when the charge level increase is proportional to the charging time, i.e., when:
\begin{equation}
    \label{eq:linear_phi}
    \varphi(\ell,\tau) = \max\big(0,\min(\mc, \ell + \alpha \tau )\big)  \quad \text{for some $\alpha > 0$.}
\end{equation}

\noindent \textbf{Optimal schedule characterization, simple formula for $\fc{a}$ and $\bc{a}$.}
Let $a$ be the station sequence $u = s_0,s_1, \ldots, s_k,s_{k+1} = v$ be a station sequence between $u$ and $v$.
We use the convention $\fu[u] = 0$ if $u \in \De$, and introduce the constants
\begin{align*}
    \lav &= \alpha \overbrace{\Big(\tb[v] - \te[u] - \displaystyle\sum_{i=1}^{k+1} \ti[s_{i-1},s_i]\Big)}^{\text{Charging time available during $a$}}  \quad \text{and}\quad
    \ldinc = \overbrace{\lav}^{\substack{\text{Maximum} \\ \text{recharge} \\ \text{during $a$}}}   \overbrace{-\displaystyle\sum_{i=1}^{k+1} \fu[s_{i-1},s_{i}] - e_v}^{\text{Energy consumed during $a$}},
\end{align*}
Constant $\ldinc$ is interesting because, if battery level was not lower bounded by $0$ and upper bounded by $\mc$, i.e., if we had $\varphi(\ell,\tau) = \ell + \alpha \tau$, we would have $\fc{a}(\lin) = \lin + \ldinc$ and $\bc{a}(\lout) = \lout - \ldinc$. 
The following proposition introduces adapted closed formula for $\fc{a}$ and $\bc{a}$ under linear recharge~\eqref{eq:linear_phi}.
}

\myblue{
\begin{theorem}
\label{thm:stationSequenceFu}
Let $\ldi$ be defined as $\fc{a}\big(\lmb\big) - \lmb$.
We have
\begin{alignat*}{1}
\fc{a}(\lin) &= 
\begin{cases}
\min\big(\lme,\lin + \ldi\big) & \text{if } \lin \geq \lmb, \\
-\infty & \text{otherwise,}
\end{cases}
\\
\bc{a}(\lout) &= 
\begin{cases}
\lmb + \max(0, \lout - \ldi) & \text{if } \lout \leq \lme, \\
+\infty & \text{otherwise,} 
\end{cases}
\end{alignat*}
and the following closed formulas for   $\lmb$, $\lme$, and $\ldi$ 
\begin{alignat*}{1}
\lmb &= \max\big(\fu[u,s_1],-\ldinc) 
 \\
\lme &=
\min\bigg(\mc, \mc+\lav - \displaystyle\sum_{i=1}^{k}\fu[s_{i-1},s_i]\bigg) - \fu[s_k,v] - e_v, \\
\ldi &= \min\big(\ldinc, \lme -\lmb \big).
\end{alignat*}
\end{theorem}

Theorem~\ref{thm:stationSequenceFu} makes battery level computation quicker: We can compute $\lmb$, $\lme$, and $\ldi$  during the preprocessing and then have simple closed formula for $\fc{a}$ and $\bc{a}$.
Note that $\lmb$ is the \emph{minimum initial charge $\ib$ such that $\ro$ is feasible given $\ib$},
 $\lme$ is the \emph{charge at the end of $a$ if the battery is fully charged at the beginning}, and $\ldi$ is the \emph{difference between the charge at the end $a$ and the charge at the beginning if the vehicle is initially at the minimum charge}~$\lmb$.
 
\vspace*{1em}
\noindent \textbf{Station sequences dominance.}
The following corollary shows that, when recharge is linear, we can  easily check dominance of station sequences in the linear case in point 2 of Theorem~\ref{prop:stationSequenceDominance}.
\begin{corollary}
\label{cor:linearDominance}
Let $a$ and $a'$ be the two station sequences $u=s_0,s_1,\ldots,s_k, s_{k+1} = v$ and $u=s'_0,s'_1,\ldots,s'_{k'}, s'_{k'+1} = v$ between $u$ and $v$ that coincide starting from one station, i.e., such that there exists $i$ and $i'$ such that $k-i = k'-i'$ and $s_{i + j} = s'_{i' + j}$ for $j$ in $\{0,\ldots,k-i\}$. 
Under linear recharge, we have $\big(\levb[s_i][a,*],\tb[s_i][a,*](\lin)\big) 
\succeq 
\big(\levb[s'_{i'}][a',*],\tb[s'_{i'}][a,*](\lin)\big)$ for all $\lin \in [0,\mc]$ if and only if:
$$ \fu[u,s_1] \leq \fu[u,s'_1] \quad \text{and} \quad \sum_{j=1}^i \fu[s_{j-1},s_j] + \alpha \ti[s_{j-1},s_{j}] \leq  \sum_{j=1}^{i'} \fu[s_{j-1},s_j] + \alpha \ti[s_{j-1},s_{j}].$$
\end{corollary}
}

\section{\myblue{Solution Methods}}
\label{sec:solution_methods}

\myblue{Now that we are able to characterize all meaningful (recharging and non-recharging) arcs from one service location or depot to another, we can turn our attention towards the solution of Problem~\eqref{eq:masterPb}. As the set of possible routes $\fR$ grows exponentially with the number of services in~$\Tr$, enumeration is generally impractical and we must rely on a dynamic generation of relevant routes during the solution process. We therefore develop a branch-and-bound algorithm, which solves the linear relaxation of Problem~\eqref{eq:masterPb} at each node by column generation, using a pricing algorithm to detect columns (i.e., routes) of negative reduced cost. This leads to a Branch-and-Price (B\&P) algorithm \citep[see, e.g.,][]{Barnhart1998} capable of finding optimal solutions to the EVSP. In the remainder of this section, we will describe the branching strategy adopted as well as the pricing algorithms. Next, we will introduce graph sparsification and diving techniques that permit to greatly speed up the solution process, leading to a fast heuristic applicable to larger instances.}

\subsection{\myblue{Branching Strategy}}
\label{sec:branching_strategy}

As seen in previous studies \citep[see, e.g., ][]{martinelli2011branch}, it is not possible to branch on route variables~$x_r$ since a branch with $x_r = 0$ would lead to this variable being generated again by the pricing subproblem. \myblue{For this reason, given a fractional solution to the linear relaxation, the algorithm uses two different branching rules. The first branching rule calculates the number of vehicles leaving a depot and branches on the depot with the most fractional value. The second is used if there is any fractional traversal from a service $\tr_1 \in \Tr$ to service $\tr_2 \in \Tr$. Analogously, it chooses the traversal with the most fractional value. Both branching rules include a constraint in each of the two new formulations generated, the first rounding down and the second rounding up the fractional value. These constraints add new dual values, which the pricing algorithm must consider. Finally, the branching tree is explored following a best-bound strategy.}

\subsection{\myblue{Pricing Algorithm}}
\label{sec:pricing_subproblem}

\myblue{
Each time a node is considered during branch-and-bound, the linear relaxation of Problem~\eqref{eq:masterPb} is solved subject to the constraints associated with the current branch. The dual values $\lambda_v$ associated to Equation~\eqref{eq:partconstraint} serve to define one pricing subproblem for each depot $\de \in \De$, whose goal is to find a column of minimum reduced cost, i.e.,
\begin{equation}\label{eq:pricingSubproblem}
	\min_{\ro \in \fR_{\de}} c_{\ro} - \sum_{\tr \in \ro} \lambda_{\tr},
\end{equation}
where $\fR_{\de}$ is the set of feasible routes that start and end in $\de$. 
We therefore iteratively solve the set of pricing problems for all depots, add up to the best 200 columns identified for each depot in the formulation, and solve the linear relaxation again. This iterative column generation process continues until no additional column of negative reduced cost can be identified. To improve convergence, we also apply a simple dual stabilization technique \citep{pessoa2010exact}, using a factor $\alpha \in [0, 1[$ to calculate the convex combination of the dual values from the last and the current linear relaxation solutions. Factor $\alpha$ starts with $0.9$ and, if the pricing fails for every depot, it is reduced by 0.1. Stabilization is finished when $\alpha = 0.0$.
}

\myblue{
In the rest of this section, we focus on the pricing problem for a given depot $\de$ and omit the subscript $\de$ in $D_{\de} = (V_{\de},A_{\de})$.
We define the reduced cost $\tilde c_a$ and $\tilde c_P$ of a $u$-$v$ arc $a$ and of a path $P$ in $D$ as:
$$\tilde{c}_a = 
\begin{cases}
    \cv + c_a - \lambda_v, & \text{if } u = o \\
    c_a, & \text{if } v=d, \\
    c_a - \lambda_v, & \text{otherwise.}
\end{cases} \quad \text{and} \quad \tilde c_P = \sum_{a \in P} \tilde c_a.
$$
Let us denote by $\calP_{od}$ the set of $o$-$d$ paths in $D$. A path $P$ in $\calP_{od}$ is a route $r$.
With the above definition of reduced costs, we have 
$ \tilde c_P := \sum_{a \in P} \tilde c_a =  c_{\ro} - \sum_{\tr \in \ro} \lambda_{\tr}$, and 
we can reformulate the pricing subproblem as the resource-constrained shortest path problem:
\begin{equation}
\label{eq:pricingRCSP}
    \min \big\{\tilde c_P \colon P \in \calP_{od}, \quad P=( o,a_,\ldots,a_{k+1},d) , \, \text{ and } \,  \fc{a_{k+1}}\circ \ldots \circ \fc{a_1}(\mc) \geq 0\big\}.
\end{equation}
Given a path $P = a_1,\ldots,a_k$ in $\calP_{od}$ and $i \in \{0,\ldots,k\}$, Lemma~\ref{lem:bidirectionalRouteFeasibility} ensures that it is a feasible solution of~\eqref{eq:pricingRCSP} if and only if:}
\begin{equation}
    \label{eq:bidirectionalPathFeasibility}
    \myblue{
    \fc{a_i}\circ \ldots \circ \fc{a_1}(\mc) \geq \bc{a_{i+1}}\circ\ldots \circ \bc{a_{k}}(0).}
\end{equation}

\myblue{The rest of this section introduces our algorithm for the resource-constrained shortest path problem~\eqref{eq:pricingRCSP}. 
It follows a classic path enumeration strategy~\citep{Irnich2005} using dominance on forward paths and bounds on backward paths to discard partial solutions. Moreover, its originality lies in the way bounds are computed on backward paths \citep{Parmentier2018} and in the use of these bounds to select in which order the paths are enumerated.}

\vspace*{1em}
\myblue{
\noindent 
\textbf{Forward resources and dominance.}
Let us define the \emph{forward resource set} \mbox{$\fQ = \bbR \times [0,M]\cup \{\fone\}$}. 
Resource $\fone$ enables to encode paths infeasibility.
The \emph{forward resource}, or simply the \emph{resource} $\fq_P$ of an $o$-$v$ path $P=a_1,\ldots,a_i$ is:
$$\fq_P = \begin{cases}
    (\tilde c_P,l_P), & \text{if } l_P \geq 0 \\
    \fone, & \text{otherwise,} 
\end{cases} 
\quad \text{where} \quad
\tilde c_P = \sum_{a \in P}\tilde c_a 
\quad\text{and}\quad
l_P = \fc{a_i}\circ \ldots \circ \fc{a_1}(\mc). $$
For each arc $a$, we introduce the forward resource extension function $\fref_a : \fQ \rightarrow \fQ$ as:
$$ \fref_a(\fq) = \begin{cases}
    (\tilde c + \tilde c_a, \fc{a}(l)), 
    & \text{if} \quad \fq = (\tilde c,l) 
    \quad\text{and}\quad \fc{a}(l) \geq 0, \\
    \fone & \text{otherwise.}
\end{cases} $$
Given an $o $-$v$ path $P = Q+a$ composed of path $Q$ followed by arc $a$, it follows from these definitions that $\fq_P = \fref_a(\fq_Q)$. We introduce the partial order $\preceq$ on $\fQ$ defined by $\fq \preceq \fone$ for all $\fq$ in $\fQ$, and:
$$(\tilde c, l) \preceq (\tilde c',l') \quad \text{if and only if} \quad c \leq c' \quad \text{and} \quad l \geq l'.$$
\myred{Note that this order is unrelated to the one introduced in Section~\ref{sec:StationSequenceScheduling} (these relations apply to different resource sets, so their meaning is clear from context)}.
It is clear that forward resource extension functions $\fref_{a}$ are non-decreasing with respect to $\preceq$.
We say that an $o$-$v$ path $P$ \emph{dominates} an $o$-$v$ path $P'$ if $\fq_P \preceq \fq_{P'}$. The following lemma shows that we can safely ignore dominated paths.
}
\myblue{
\begin{lemma}\label{lem:dominance}
    Suppose that an optimal solution $P = Q + R$ of~\eqref{eq:pricingRCSP} is composed of an $o$-$v$ path $Q$ followed by a $v$-$d$ path $R$. If $Q'$ is an $o$-$v$ path that dominates $Q$, then the $o$-$v$ path $P' = Q' + R$ composed of $P'$ followed by $R$ is an optimal solution of~\eqref{eq:pricingRCSP}. 
\end{lemma}
}

\myblue{
\textsc{Proof.}
Suppose that $Q$, $Q'$, and $R$ are of the form $o,a_1,\ldots,a_i,v$, \quad $o,a'_1,\ldots,a'_{i'},v$, and $v,a_{i+1},\ldots, a_{k+1},d$, respectively.
Since $P$ is feasible $\fc{a}$ is non-decreasing for all $a$, we have 
$0\leq \fc{a_{k+1}} \circ \fc{a_1}(\mc) 
= \fc{a_{k+1}} \circ \fc{a_{i+1}}\big(\fc{a_{i}} \circ \fc{a_1}(\mc)\big) 
\leq \fc{a_{k+1}} \circ \fc{a_{i+1}}\big(\fc{a'_{i'}} \circ \fc{a'_1}(\mc)\big)$, 
and $P'$ is feasible.
Besides $\tilde{c}_{P'} = \tilde{c}_{Q'} + \tilde c_R \leq \tilde{c}_{Q} + \tilde c_R  = \tilde c_P$, which gives the result.\halmos}

\vspace*{1em}
\myblue{
\noindent \textbf{Backward resources and lower bounds.}
Let us define the \emph{backward resource set} $\bQ = \bbR \times [0,M]\cup \{\bone\}$. 
The \emph{backward resource}, or simply the \emph{resource} $\bq_P$ of a $v$-$d$ path $P=a_1,\ldots,a_i$~is:
$$\bq_P = \begin{cases}
    (\tilde c_P,m_P), & \text{if } m_P \leq M \\
    \bone, & \text{otherwise,} 
\end{cases} 
\quad \text{where} \quad
\tilde c_P = \sum_{a \in P}\tilde c_a 
\quad\text{and}\quad
m_P = \bc{a_{1}}\circ\ldots \circ \bc{a_{i}}(0).$$
}\myblue{Again, $\bone$ encodes paths infeasibility.
The component $m_P$ can be interpreted as the minimum battery level required at the beginning of $P$ to operate $P$ without running out of battery.
For each arc $a$, we introduce the backward resource extension function $\bref_a : \bQ \rightarrow \bQ$ as:
$$ \bref_a(\bq) = \begin{cases}
    \bone, & \text{if} \quad \bq = \fone 
    \quad \text{or}\quad \bq = (\tilde c,m) 
    \quad \text{and}\quad \bc{a}(m)) = +\infty \\
    (\tilde c + \tilde c_a, \bc{a}(m)), 
    & \text{if} \quad \bq = (c,\ell) 
    \quad\text{and}\quad \bc{a}(m) \leq M.
\end{cases}$$
}\myblue{Given a $v$-$d$ path $P = a+R$ composed of and arc $a$ followed by a path $R$, it follows from these definitions that $\bq_R = \bref_a(\bq_Q)$. 
We introduce the partial order $\preceq$ on $\bQ$ defined by $\bq \preceq \bone$ for all $\bq$ in $\bQ$, and:
$$(\tilde c, m) \preceq (\tilde c',m') \quad \text{if and only if} \quad \tilde c \leq \tilde c' \quad \text{and} \quad m \leq m'. $$
$(\bQ,\preceq)$ happens to be a \emph{meet semi-lattice}: any pair $(\bq,\bq')$ of elements of $\bQ$ admits a greatest lower bound or \emph{meet} $\bq \meet \bq'$. We have $\bq \meet \bq' = \bq'$ if $\bq = \bone$, and $\bq \meet \bq' = \bq$ if $\bq' = \bone$, and:
$$ (\tilde c, m) \meet (\tilde c, m') = \big(\min(\tilde c, \tilde c'), \min (m,m')\big).$$
}\myblue{We can now introduce the following generalization of the dynamic programming equation:
\begin{equation}\label{eq:generalizedDynamicProgramming}
    b_v = \begin{cases}
        (0,0), & \text{if } v=d, \\
        \displaystyle\bigmeet_{a = (v,w) \in \delta^+(v)} \bref_a(b_w), & \text{otherwise}
    \end{cases}
\end{equation}
where the meet operator $\meet$ plays the role played by the minimum in classic dynamic programming algorithm. 
Using the classic dynamic programming algorithm on acyclic digraphs,
its unique solution can be computed recursively along a reverse topological order.
An induction along a reverse topological order gives the following result (\citealt{Parmentier2018}, Theorem~3).}

\myblue{
\begin{lemma}\label{lem:bound}
    The unique solution of Equation~\eqref{eq:generalizedDynamicProgramming} is such that:
    \begin{equation}\label{eq:boundProperty}
        b_d = (0,0) \quad \text{and} \quad b_v \preceq \bq_R \quad \text{for any vertex $v$ and $v$-$d$ path $R$}.
    \end{equation}
\end{lemma}
}

\myblue{
\begin{remark}
    A classic labeling approach for this problem would consist in computing separately lower bounds for each resource using a (traditional) dynamic programming approach~\citep{Dumitrescu2003}. 
    In the case of the EVSP, this would amount to running a separate shortest path algorithm for the reduced costs and the battery level.
    The bounds obtained with our approach are better because whenever we identify that the battery is empty, we have $\bref_a(b_w) = \bone$, and the corresponding reduced cost is not taken into account in the~$\meet$ operator, whereas it is taken into account in the classic approach.\\
\end{remark}
}

\myblue{
\begin{remark}
    With the classic dynamic programming equation, the existence of at least one $\argmin$ in the minimum enables to rebuild the solution. 
    As there is no notion of $\argmin$ for the meet operator $\meet$,  Equation~\eqref{eq:generalizedDynamicProgramming} permits us to compute a bound but not an optimal solution.
\end{remark}
}

\vspace*{1em}
\myblue{
\noindent 
\textbf{Bidirectional feasibility check.}
Let us finally define the bidirectional resource function~$\fbc : \fQ \times \bQ \rightarrow \bbR \cup \{+ \infty\}$ as:
$$ \fbc(\fq,\bq) = \begin{cases}
    \fvec{c} + \bvec{c} & \text{if} \quad \fq = (\fvec{c},l), \quad \bq = (\bvec{c},m), \quad \text{and} \quad l \geq m, \\
    +\infty & \text{otherwise.}
\end{cases} $$
As a consequence of Lemmas~\ref{lem:bidirectionalRouteFeasibility} and \ref{lem:bound}, we obtain the following result.
}

\myblue{
\begin{lemma}
\label{lem:bidirectionalPath}
    Let $Q$ be an $o$-$v$ path and $b_v$ a lower bound satisfying Equation~\eqref{eq:boundProperty}. Then $\fbc(q_Q,b_v) \leq \tilde c_P$ for any feasible $o$-$d$ path $P = Q +R$ starting by $P$.
\end{lemma}}

\myblue{
\textsc{Proof.} Let $\fq{Q} = (\tilde c_Q, l_Q)$ and $\bq{R} = (\tilde c_R, l_R)$ be the resources of $Q$ and $R$, and $b_v = (\tilde c_v,l_v)$.
Since $P$ is feasible, we have $l_Q \geq l_R \geq l_v$. 
Hence $\fbc(q_Q,b_v) = \tilde c_Q + \tilde c_v \leq \tilde c_Q + \tilde c_R = \tilde c_P$, which gives the result.
\halmos
}

\myblue{
Hence, if we already know a feasible $o$-$d$ path $P$ with reduced cost $\tilde c$, and $\fbc(q_Q,b_v)\geq \tilde c$, then we know that no path starting by $Q$ is going to be better $P$, and we can safely discard $Q$.
}

\vspace*{1em}
\myblue{
\noindent \textbf{Algorithm.}
Assuming that we have computed lower bounds $b_v$ for $v$ in $V$ satisfying  Equation~\eqref{eq:boundProperty} during a preprocessing phase, we propose Algorithm~\ref{alg:enumeration} for the pricing problem. It implicitly enumerates all the $o$-$v$ paths that cannot be discarded using dominance (Lemma~\ref{lem:dominance}) or bounds (Lemma~\ref{lem:bidirectionalPath}).
The algorithm maintains a heap $\sfL$ containing $o$-$v$ paths $P$ for any $v$, a list $\mathsf{L}_{v}^{\mathrm{nd}}$ containing $o$-$v$ paths for each vertex $v$ in $V\backslash \{o,d\}$, the best solution found $\mathtt{current\_best\_solution}$ and its cost $\mathtt{current\_best\_cost}$. 
Any $o$-$v$ path $P$ in $\sfL$ is stored with its resource $\fq_P$ and its \emph{key} $k_P = \fbc(\fq_Q,b_v)$. By Lemma~\ref{lem:bidirectionalPath}, $k_P$ is a lower bound on the cost of any $o$-$d$ path starting by $P$. 
And since $b_d =(0,0)$, if $P$ is in an $o$-$d$ path, $k_P = (\tilde c_P)$.
}

\begin{algorithm}[htbp]
\linespread{1}
\myblue{
    \begin{algorithmic}[1]
    \STATE \textbf{input}: an acyclic digraph $D=(V,A)$, resource extension functions $\fref_a$, bounds $(b_v)_{v \in V}$ satisfying~\eqref{eq:boundProperty}; 
    \STATE \textbf{output}: a path with minimum reduced cost
    \STATE \textbf{initialization:}  $\mathtt{current\_best\_cost} \leftarrow +\infty$, $\mathtt{current\_best\_solution} \leftarrow \mathtt{undefined}$, $\mathsf{L} \leftarrow \emptyset$, and $\mathsf{L}_{v}^{\mathrm{nd}} \leftarrow \emptyset$ for each $v \in V \backslash \{o,d\}$;
    \STATE $P_o \leftarrow $ empty path  at the origin $o$;
    \STATE $\fq_{P_o}\leftarrow (0,M)$; \quad $k_{P_o}\leftarrow\fbc\big(\fq_o,b_o\big)$;
    \STATE add $P_o$ to $\mathsf{L}$ with resource $\fq_o$ and key $k_{P_o}$;
    \WHILE{$\mathsf{L}$ is not empty}
    \STATE $P \leftarrow$ a path $\tilde P$ of minimum key $k_{\tilde P}$ in $\mathsf{L}$; \label{step:pathSelection} 
    \STATE remove $P$ from $\sfL$;
    \STATE $v\leftarrow$ last vertex of $P$; \label{step:pathSelectionEnd}
    \FORALL{$a\in\delta^+(v)$}\label{step:extensionStarts}
    \STATE $Q \leftarrow P+a$; \label{step:considerStarts}
    \STATE $w\leftarrow $ last vertex of $Q$; \label{step:Qstarts}
    \STATE $\fq_Q \leftarrow \fref_a(\fq_P)$; \quad $k_{Q} \leftarrow \fbc(\fq_Q,b_w)$;
    \IF{$k_Q < \mathtt{current\_best\_cost} $} \label{step:LBtest}
    \IF{$w=d$} \label{step:extStart}
    \STATE $\mathtt{current\_best\_cost} \leftarrow k_Q $, \quad  $\mathtt{current\_best\_solution} \leftarrow Q $; \label{step:codub}
    \STATE remove from $\sfL$ all the paths $R$ with $k_{R} \geq \mathtt{current\_best\_cost}$
    \ELSE  
    \IF{$Q$ is not dominated by any path in $\mathsf{L}_{w}^{\mathrm{nd}}$} \label{step:DomTest} 
    \STATE $\mathsf{L}_{w}^{\mathrm{nd}} \leftarrow \mathsf{L}_{w}^{\mathrm{nd}} \cup \{Q\}$ and remove from $\mathsf{L}_{w}^{\mathrm{nd}}$ and $\mathsf{L}$ every path dominated by $Q$; \label{step:removal}
    \STATE add $Q$ to $\mathsf{L}$ with resource $\fq_Q$ and key $\fbc(\fq_Q,b_w)$; \label{step:updateL}
    \ENDIF 
    \ENDIF \label{step:Qends}
    \ENDIF \label{step:considerEnds}
    \ENDFOR \label{step:extensionEnds}
    \ENDWHILE
    \STATE \textbf{return} $\mathtt{current\_best\_solution}$;
    \end{algorithmic}
    \caption{Enumeration algorithm for the pricing problem~\eqref{eq:pricingRCSP}}
    \label{alg:enumeration}
}
\end{algorithm} 

\myblue{
At a given iteration, the algorithm extracts from $\sfL$ an $o$-$v$ path $P$ with minimum key $k_P = \fbc(\fq_P,b_v)$ (steps~\ref{step:pathSelection}-\ref{step:pathSelectionEnd}), and \emph{extends} it (steps~\ref{step:extensionStarts}-\ref{step:extensionEnds}). 
That is, for each arc $a$ outgoing from~$v$, the path $Q=p+a$ is \emph{considered} by the algorithm (steps~\ref{step:considerStarts}-\ref{step:considerEnds}). 
We start by checking if $Q$ can be safely discarded because $k_Q \geq \mathtt{current\_best\_cost}$ (Lemma~\ref{lem:bidirectionalPath}).
If $Q$ is not discarded and is an $o$-$d$ path, since $k_Q = \tilde c_Q$, it means that it is the best solution found, and we update $\mathtt{current\_best\_solution}$ and $\mathtt{current\_best\_cost}$ accordingly.
Otherwise, we check if $Q$ can be discarded because it is dominated by a previously considered $o$-$v$ path in $\mathsf{L}_{v}^{\mathrm{nd}}$ (Lemma~\ref{lem:dominance}). If not, we update $\mathsf{L}_{v}^{\mathrm{nd}}$ with $Q$ and add it to the heap $\sfL$ of paths to be extended.
}

\myblue{
At any given time in the algorithm, $\sfL$ contains the paths which have been considered without being discarded by the algorithm but that have still not been extended. For each vertex~$v$, the list $\sfL_v$ contains the $o$-$v$ paths which have been considered by the algorithm and are not dominated by any other $o$-$v$ path considered by the algorithm.
}

\myblue{
\begin{theorem}\label{theo:AlgoConvergence}
    Algorithm~\ref{alg:enumeration} converges after a finite number of iterations. At the end of the algorithm, $\mathtt{current\_best\_solution}$ contains an optimal solution of Equation~\eqref{eq:pricingRCSP} if it admits a feasible solution, and is equal to $\mathtt{undefined}$ otherwise.
\end{theorem}
}

\myblue{
The proof of Theorem~\ref{theo:AlgoConvergence} (\citealt{Parmentier2018}, Theorem 1) relies on the following invariant. Let~$P$ be a feasible $o$-$d$ path whose subpaths are not dominated. Then at any time during the algorithm, either $P$ or one of its subpaths is in $\sfL$, or $\mathtt{current\_best\_solution}$ contains a feasible path $Q$ whose reduced cost $\tilde c_Q$ is non-greater than $\tilde c_P$.
}

\subsection{Graph Sparsification}

\myblue{
As seen in our experiments, the B\&P algorithm described earlier in this section can find optimal solutions for instances with 100 services but does not scale up to larger cases with 500 services in a reasonable time. To speed up our approach in such cases, we introduce graph sparsification strategies that permit to reduce the number of arcs in the digraphs $D_{\de} = (V_{\de},A_{\de})$ for each depot $\de$. This simplification has only a moderate impact on solution quality while permitting to solve much larger problems.
}

\myblue{
As previously, to make the notations lighter, we drop the depot in subscript and use the notations $D=(V,A)$. We use three main intuitions to select arcs. (i) An arc $a = (u,v)$ is interesting if it does not require driving an unnecessary distance ($c_a$ is small), and does not make the vehicle waste time waiting for a service ($\tb[v] - \te[u]$ is small). 
(ii)~On arcs with recharge, some time must be spent recharging, but it is useless to wait at the station without recharging.
 To measure this, given an arc $a$ corresponding to a station sequence $u=s_0,s_1,\ldots,s_{k},s_{k+1}=v$ between $u$ and $v$,
 we define the \emph{lost time} $L_a$ as the duration between the time when the vehicle reaches the maximum charge in station $s_k$ and the time when it leaves $s_k$  under an optimal scheduling if the vehicle starts $a$ with the minimum charge level $e_{s_0,s_1}$ required to reach the first station.
    When recharge is linear, this amount can be computed with the formula:
    $$ L_a = \max \bigg( 0, \tb[s_1] - \Big(\te[s_1] + \underbrace{\sum_{i=s}^{k+1} \ti[s_{i-1},s_i] }_{\text{travel time}} + \underbrace{\frac{\mc + \sum_{i=2}^{k}\fu[s_{i-1},s_i]}{\cf}}_{\text{maximum charging duration}} \Big)\bigg).$$
Finally, (iii)~some arcs from and to the depot must be maintained to ensure the feasibility of the instances.

To construct the sparsified set of arcs, we partition the original arcs $A$ into four subsets $A_1 \cup A_2\cup A_3 \cup A_4$, each subset of arcs having different properties.
For each subset $A_i$, we define a goodness measure $\gamma_i : A_i \rightarrow \bbR$.
We then build a sparsified version $A'_i$ of $A_i$ as follows:
$$ A'_i = \left\{a = (u,v) \in A_i \colon
\begin{array}{l}
\text{$a$ is among the $\nu_i$ arcs with smallest $\gamma_i(a)$ in $\delta^+(u) \cap A_i$ }\\
\text{ or in the  $\nu_i$ arcs with smallest $\gamma_i(a)$ in $\delta^-(v) \cap A_i$}
\end{array}
\right\}.$$
Ties are broken randomly in the presence of multiple arcs with identical $\gamma_i(a)$ values. Finally, we define the sparsified subset as $A' = A'_1 \cup A'_2\cup A'_3 \cup A'_4$.
The definition of $A'_i$ ensures that sufficiently many good incoming arcs and outgoing arcs of $A_i$ are preserved for each vertex, therefore giving enough flexibility in the route choices. 

We recall that an arc $a$ corresponds to a station sequence, with potentially no station. Let $A^{\mathrm{s}}$ denote the set of arcs whose station sequence contains at least one station in $\St$. Let us now define the $A_i$, $\gamma_i$ and $\nu_i$ as follows:
\begin{align*}
    A_1 &= \big\{a=(u,v) \in A \backslash A^{\mathrm{s}} \colon u \in \De \text{ or } v\in \De \big\} &  
    \nu_1&  = 2 &  
    \gamma_1(a) &= \tilde c_a \\
    A_2 &= \big\{a=(u,v) \in A^{\mathrm{s}} \colon u \in \De \text{ or } v\in \De \big\} &  
    \nu_2&  = 2 &  
    \gamma_2(a) &= \tilde c_a \\    
    A_3 &= \big\{a=(u,v) \in A \backslash A^{\mathrm{s}} \colon u \in \Tr \text{ and } v\in \Tr \big\} &  
    \nu_3&  = 15 &  
    \gamma_3(a) &= \tilde c_a + 0.1 (\te[v] - \tb[u]) \\ 
     A_4 &= \big\{a=(u,v) \in A^{\mathrm{s}} \colon u \in \Tr \text{ and } v\in \Tr \big\} &  
    \nu_4&  = 2 &  
    \gamma_4(a) &= \tilde c_a + 0.1 (\tb[v] - \te[u]) + 0.1(L_a).
\end{align*}
The values of the different constants have been selected to achieve a good balance between the different factors. Notably, the larger value associated with $\nu_i$ reflects the fact that most arcs in the optimal solution are going to be arcs without recharge in $A_3$.
}

\subsection{Diving Heuristic}
\label{sec:diving}

Diving heuristics have been successfully applied to many complex combinatorial optimization problems \citep{sadykov2019primal}, \myblue{where they permit to quickly find good primal solutions.} In its most basic version, a diving heuristic consists of a depth-first exploration of the B\&P tree. At each search node, the column generation is first completed, and the route variable $x_r$ with the largest value is fixed to $1$ (ties are broken randomly). The diving algorithm terminates whenever an integer solution is achieved or the column generation does not return any feasible solution.

To speed up the solution process, we use the sparsification strategy described in the previous section, and adapt the diving approach to fix multiple routes at each search node, prioritizing those with variables $x_r$ closest~to~$1$. To that end, we consider the routes in non-decreasing order of their associated $x_r$ values. Each route according to this order is fixed if it does not contain a service from a previously fixed route, otherwise it is skipped. This iterative process is pursued until $n_\textsc{dmr}$ routes have been fixed for this node or all routes have been considered.

The resulting diving algorithm is fast but greedy since integer columns are permanently fixed. To improve solution quality, we rely on a more sophisticated approach called strong diving \citep{sadykov2019primal}. At each search node, this algorithm variant evaluates several possible candidate routes instead of fixing a single one. For each such candidate, it tentatively fixes the associated route and calculates the result of the column generation. The candidate leading to the smallest linear relaxation after fixing is selected, and the search proceeds to the next node. Preliminary experiments revealed that this approach helps drive some of the method's early decisions, but its systematic application to all search nodes can lead to computational time overhead. Therefore, in the proposed algorithm, strong diving is used only for a fixed number of steps, after which the method follows up with the regular diving heuristic.

\section{Experimental Analyses}
\label{sec:numerical_results}

This section analyses the performance of the proposed methods and the impact of key methodological choices related to the sparsification and diving algorithm. First, we present the results of the B\&P algorithm, considering the regular version and its version using sparsification. Then, we focus on the diving heuristics results and draw a sensitivity analysis of its key parameters. These tests are conducted on the \emph{large} benchmark instance set from \cite{Wen2016}, containing instances with $\{2, 4, 8\}$ depots, $\{4, 8, 16\}$ recharging stations and $\{100, 500\}$ services.
Each instance is named following the format \texttt{Dw\_Sx\_Cy\_z}, where $w$ is the number of depots, $x$ is the number of recharging stations, $y$ is the number of services, and $z$ is an index that differentiates each instance. 

An important characteristic of this instance set is that the fixed vehicle cost is set to a large value of $\cv = 10{,}000$. Consequently, vehicle use costs dominate the rest of the objective contribution related to the distance, indirectly implying a hierarchical objective seeking fleet-size minimization in priority followed by driving cost optimization. To conduct fine-grained analyses, we will therefore report the original cost of the solutions along with their fleet size and driving cost.

All algorithms have been developed in C\texttt{++}, using CPLEX 12.8 for the linear programs of the B\&P master problem. Note that in the diving algorithm, we rebuild the mathematical model after each route fixing, as this progressively reduces the size of the master formulation and improves the overall solution time. We conduct all experiments on a computer with an Intel Core i7-8700K CPU @\,3.70GHz and 64 GB of RAM, running Ubuntu Linux in a single thread. A time limit of six hours (21,600 seconds) was used for each run.

\subsection{Exact Solution with the Branch-and-Price Algorithm}
\label{sec:BPresults}

We developed two versions of the B\&P algorithm. The first solves a pricing subproblem on the complete underlying graph, whereas the second version solves the pricing problem only on the sparsified graph. The results for the B\&P algorithm on the instances with 100 and 500 services are presented in Tables \ref{tbl:bap100} and \ref{tbl:bap500}.
The first four columns of Table \ref{tbl:bap100} report for each instance the best results over five runs from the adaptive large neighborhood search (ALNS) of \citet{Wen2016}: the upper bound (UB), the number of vehicles used (V), the driving cost (D), and the average computational time (T) for each run in seconds (ALNS runs were performed on an Intel i7-3520M @\,2.9 GHz, which is $1.6\times$ slower than our processor based on Passmark single-thread CPU ratings).
The next two groups of five columns give the results of the regular B\&P and the B\&P with sparsification: the optimal solution found (Opt), the number of vehicles~(V), the driving cost~(D), the computational time~(T) in seconds, and the number of B\&B nodes (\#). The best solution value for each instance is highlighted in boldface.

\begin{table}[htbp]%
\centering
\setlength{\tabcolsep}{4pt}
\resizebox{\textwidth}{!}{%
\begin{tabular}{lccccccccccccccccccc}\toprule
\multirow{2}{*}{Instance} & \multicolumn{4}{c}{ALNS} && \multicolumn{6}{c}{Branch-and-Price} && \multicolumn{6}{c}{Branch-and-Price + Sparsification} \\ \cline{2-5}\cline{7-12}\cline{14-19}
 & \multicolumn{1}{c}{UB} & \multicolumn{1}{c}{V} & \multicolumn{1}{c}{D} &  \multicolumn{1}{c}{T} && \multicolumn{1}{c}{Root} & \multicolumn{1}{c}{Opt} & \multicolumn{1}{c}{V} & \multicolumn{1}{c}{D} & \multicolumn{1}{c}{T} & \multicolumn{1}{c}{\#} && \multicolumn{1}{c}{Root} & \multicolumn{1}{c}{Opt} & \multicolumn{1}{c}{V} & \multicolumn{1}{c}{D} & \multicolumn{1}{c}{T} & \multicolumn{1}{c}{\#} \\\midrule
\texttt{D2\_S4\_C100\_1} & 211775.0 & 21 & 1775.0 & 252 && 211734.4 & \textbf{211741.0} & 21 & 1741.0 & 22.1 & 85 && 211767.4 & 211773.9 & 21 & 1773.9 & 11.3 & 107 \\
\texttt{D2\_S4\_C100\_2} & 182178.0 & 18 & 2178.0 & 890 && 181921.4 & \textbf{181932.1} & 18 & 1932.1 & 39.2 & 217 && 181970.7 & 181987.4 & 18 & 1987.4 & 53.4 & 693 \\
\texttt{D2\_S4\_C100\_3} & 192230.0 & 19 & 2230.0 & 998 && 182227.1 & \textbf{182231.7} & 18 & 2231.7 & 8.2 & 7 && 182326.1 & 182332.4 & 18 & 2332.4 & 3.0 & 21 \\
\texttt{D2\_S4\_C100\_4} & 212231.0 & 21 & 2231.0 & 319 && 212105.7 & \textbf{212115.7} & 21 & 2115.7 & 96.5 & 109 && 212171.1 & 212175.0 & 21 & 2175.0 & 7.7 & 17 \\
\texttt{D2\_S4\_C100\_5} & 181882.0 & 18 & 1882.0 & 381 && 181679.7 & \textbf{181685.2} & 18 & 1685.2 & 12.9 & 13 && 181680.2 & \textbf{181685.2} & 18 & 1685.2 & 4.3 & 17 \\
\texttt{D4\_S8\_C100\_1} & 191600.0 & 19 & 1600.0 & 281 && 191465.5 & \textbf{191470.7} & 19 & 1470.7 & 21.6 & 31 && 191478.4 & 191492.9 & 19 & 1492.9 & 58.3 & 513 \\
\texttt{D4\_S8\_C100\_2} & 192097.0 & 19 & 2097.0 & 500 && 191897.6 & \textbf{191902.5} & 19 & 1902.5 & 22.5 & 45 && 191964.9 & 191976.1 & 19 & 1976.1 & 69.4 & 773 \\
\texttt{D4\_S8\_C100\_3} & 191510.0 & 19 & 1510.0 & 325 && 191391.1 & \textbf{191401.7} & 19 & 1401.7 & 153.9 & 399 && 191399.0 & 191406.3 & 19 & 1406.3 & 14.7 & 109 \\
\texttt{D4\_S8\_C100\_4} & 211612.0 & 21 & 1612.0 & 306 && 211461.4 & \textbf{211468.4} & 21 & 1468.4 & 50.9 & 85 && 211517.4 & 211527.3 & 21 & 1527.3 & 22.2 & 203 \\
\texttt{D4\_S8\_C100\_5} & 191704.0 & 19 & 1704.0 & 311 && 191586.8 & \textbf{191592.5} & 19 & 1592.5 & 40.1 & 89 && 191612.3 & 191614.3 & 19 & 1614.3 & 3.2 & 9 \\\midrule
Average & 195881.9 & 19.4 & 1881.9 & 456.3 && 194747.1 & 194754.1 & 19.3 & 1754.1 & 46.8 & 108.0 && 194788.8 & 194797.1 & 19.3 & 1797.1 & 24.7 & 246.2 \\\bottomrule
\end{tabular}
}
\caption{Branch-and-Price results for 100-services instance set}
\label{tbl:bap100}
\end{table}

As visible in these experiments, the proposed B\&P algorithm can solve all 100-services
instances to proven optimality within less than a minute on average. 
\myred{These instances have on average 3.34 non-dominated station sequences per pair of services $(u,v)$ such that $v$ can be operated after $u$, which leads to 14.8k arcs on average per instance}.
The optimal solutions have $0.53\%$ fewer vehicles and $6.87\%$ less driving cost than the solutions produced by the ALNS on average. For instance \texttt{D2\_S4\_C100\_3} in particular, our approach found an optimal solution with one less vehicle for a marginal increase of distance ($0.02\%$).
For these instances, the sparsification strategy permitted to reduce the computational time by $47\%$, but it also led to a $2.4\%$ increase in the driving cost (the estimated solution deterioration amounts to $0.02\%$ according to the original objective). For instance \texttt{D2\_S4\_C100\_5}, the B\&P with sparsification produced the same result as the algorithm without sparsification, in only one third of the computational time. This shows that the proposed sparsification strategy retains useful arcs and can significantly improve scalability.\\

We now focus our analyses on the large instances with 500 services. Due to their size, these instances are very challenging and cannot currently be solved to optimality within the 6 hours time limit. We will therefore focus the analysis of Table~\ref{tbl:bap500} on the lower bounds produced by the B\&P algorithm. 
Thus, the first two groups of four columns present: the root node value (Root), its time in seconds (T), the lower bound or sparsification bound at the end of the time limit (LB or SB), and the number of B\&B nodes (\#). The last columns present the percentage difference between the results of the regular B\&P and those of the B\&P with sparsification.

\begin{table}[htbp]%
\centering
\setlength{\tabcolsep}{4.5pt}
\resizebox{\textwidth}{!}{%
\begin{tabular}{lccccccccccccc}\toprule
\multirow{2}{*}{Instance} & \multicolumn{4}{c}{Branch-and-Price} && \multicolumn{4}{c}{Branch-and-Price + Sparsification} && \multicolumn{3}{c}{Difference(\%)} \\ \cline{2-5}\cline{7-10}\cline{12-14}
 & \multicolumn{1}{c}{Root} & \multicolumn{1}{c}{T} & \multicolumn{1}{c}{LB} & \multicolumn{1}{c}{\#} && \multicolumn{1}{c}{Root} & \multicolumn{1}{c}{T} & \multicolumn{1}{c}{SB} & \multicolumn{1}{c}{\#} && \multicolumn{1}{c}{Root} & \multicolumn{1}{c}{T} & \multicolumn{1}{c}{\#} \\\midrule
\texttt{D4\_S8\_C500\_1} & 835216.6 & 804.0 & 835222.8 & 1120 && 835360.0 & 63.8 & 835371.5 & 7573 && 0.02 & -92.06 & 576.16 \\
\texttt{D4\_S8\_C500\_2} & 925494.6 & 693.8 & 925500.6 & 2099 && 925875.2 & 38.2 & 925884.2 & 10735 && 0.04 & -94.50 & 411.43 \\
\texttt{D4\_S8\_C500\_3} & -- & -- & -- & -- && 766005.9 & 513.1 & 766015.2 & 1347 && -- & -- & -- \\
\texttt{D4\_S8\_C500\_4} & 825642.1 & 4041.4 & 825646.6 & 220 && 825946.4 & 1137.6 & 825959.7 & 4365 && 0.04 & -71.85 & 1884.09 \\
\texttt{D4\_S8\_C500\_5} & 786318.8 & 20564.5 & 786318.8 & 2 && 786991.2 & 785.4 & 787006.6 & 919 && 0.09 & -96.18 & 45850.00 \\
\texttt{D8\_S16\_C500\_1} & 814904.4 & 1909.8 & 814907.3 & 370 && 815116.8 & 58.1 & 815124.0 & 2557 && 0.03 & -96.96 & 591.08 \\
\texttt{D8\_S16\_C500\_2} & 824935.8 & 1178.9 & 824939.6 & 557 && 825249.7 & 85.1 & 825256.1 & 2907 && 0.04 & -92.78 & 421.90 \\
\texttt{D8\_S16\_C500\_3} & 874667.2 & 881.8 & 874670.3 & 1309 && 874958.3 & 42.9 & 874965.8 & 3882 && 0.03 & -95.13 & 196.56 \\
\texttt{D8\_S16\_C500\_4} & 784536.5 & 818.6 & 784539.5 & 1220 && 784655.6 & 41.0 & 784661.6 & 7747 && 0.02 & -94.99 & 535.00 \\
\texttt{D8\_S16\_C500\_5} & -- & -- & -- & -- && 834995.2 & 95.1 & 835003.4 & 6154 && -- & -- & -- \\\midrule
Average$^\dagger$ & 833964.5 & 3861.6 & 833968.2 & 862.1 && 834269.2 & 281.5 & 834278.7 & 5085.6 && 0.04 & -91.81 & 6308.28 \\\bottomrule
\multicolumn{14}{l}{$^\dagger$ Considering the instances solved by both approaches.}
\end{tabular}
}
\caption{Branch-and-Price results for 500-services instance set}
\label{tbl:bap500}
\end{table}

For two instances out of ten, the standard B\&P algorithm could not complete the root node's solution within the time limit. In contrast, the approach with sparsification concluded the root node for all instances. \myblue{This improvement comes from a pricing time reduction, as this component of the algorithm represents the main bottleneck of the algorithm, with 98.9\% of the computational effort dedicated to it in the approach without sparsification. In contrast, the method based on sparsification spends 82.0\% of the time on average in this procedure, and could progress much further in the column generation.}
It is noteworthy that the bound values found with sparsification are within $0.04\%$ of the original lower bounds and took $13.7\times$ less computational effort on average (considering the subset of instances for which the original algorithm could provide a bound). 
\myblue{The generation of charging arcs took less than 0.3 seconds for all instances on both approaches. It led to \myred{4.41 non-dominated station sequences per pair of services $(u,v)$ such that $u$ can be done after $v$, which gave} 452k arcs on average for the regular B\&P, and only 16k arcs on average when using sparsification.} Whereas the regular B\&P took, on average, more than one hour to solve the root node, the B\&P with sparsification took less than five minutes. Still, despite this speedup, we observed that the progress of the bounds after the root node was too slow to expect an optimal solution in a reasonable time.

\myred{Finally, given that all 100-requests instances were solved to optimality but that none of the 500-requests could be solved in a reasonable time, we conducted a final experiment with the B\&P algorithm to identify a size limit for a possible optimal solution. To do so, we extracted from each instance with 500 services three additional instances with 200, 300, and 400 random services. Considering these instances, we observed that only half of the instances with 200 services and none of the larger instances could be solved optimally in the 6-hours time limit. For larger instances beyond 100 requests, we therefore focus on diving heuristics, which will permit us to produce high-quality upper bounds.
}

\subsection{Performance of the Diving Heuristic}

The proposed diving heuristics have some parameters that need careful calibration. The ``strong diving depth limit'' (\texttt{SDDL}) parameter defines the maximum depth up to which the heuristic uses strong diving. The ``strong diving maximum solutions'' (\texttt{SDMS}) parameter limits the maximum number of solutions that are inspected during strong branching. The ``diving maximum routes'' (\texttt{DMR}) parameter sets the number of routes that are fixed in a given iteration of regular diving. Finally, the boolean parameter ``singleton routes'' (\texttt{SR}) controls the pre-generation of singleton routes, i.e., routes with only one customer, before the column generation starts. We noticed improvements with this strategy in our preliminary tests since it effectively avoids infeasible solutions after fixing variables.

Our calibration experiment was conducted on the instances with 500 services. We explored small integer values for the parameters  since our preliminary tests demonstrated that larger values led to excessive computational time amounts. We therefore consider \texttt{SDDL} $\in \{0, 1, 2, 3\}$, \texttt{SDMS} $\in \{0, 2, 3\}$, \texttt{DMR} $\in \{1, 2, 3\}$, \texttt{SR} $\in \{\text{\texttt{F}}, \text{\texttt{T}}\}$. \texttt{SDDL} and \texttt{SDMS} $ = 0$ means turning strong diving off (therefore using only regular diving). Each configuration is labeled as a quadruplet \texttt{(SDDL, SDMS, DMR, SR)}.

\begin{figure}[htbp]
    \centering
    \includegraphics[width=0.8\textwidth]{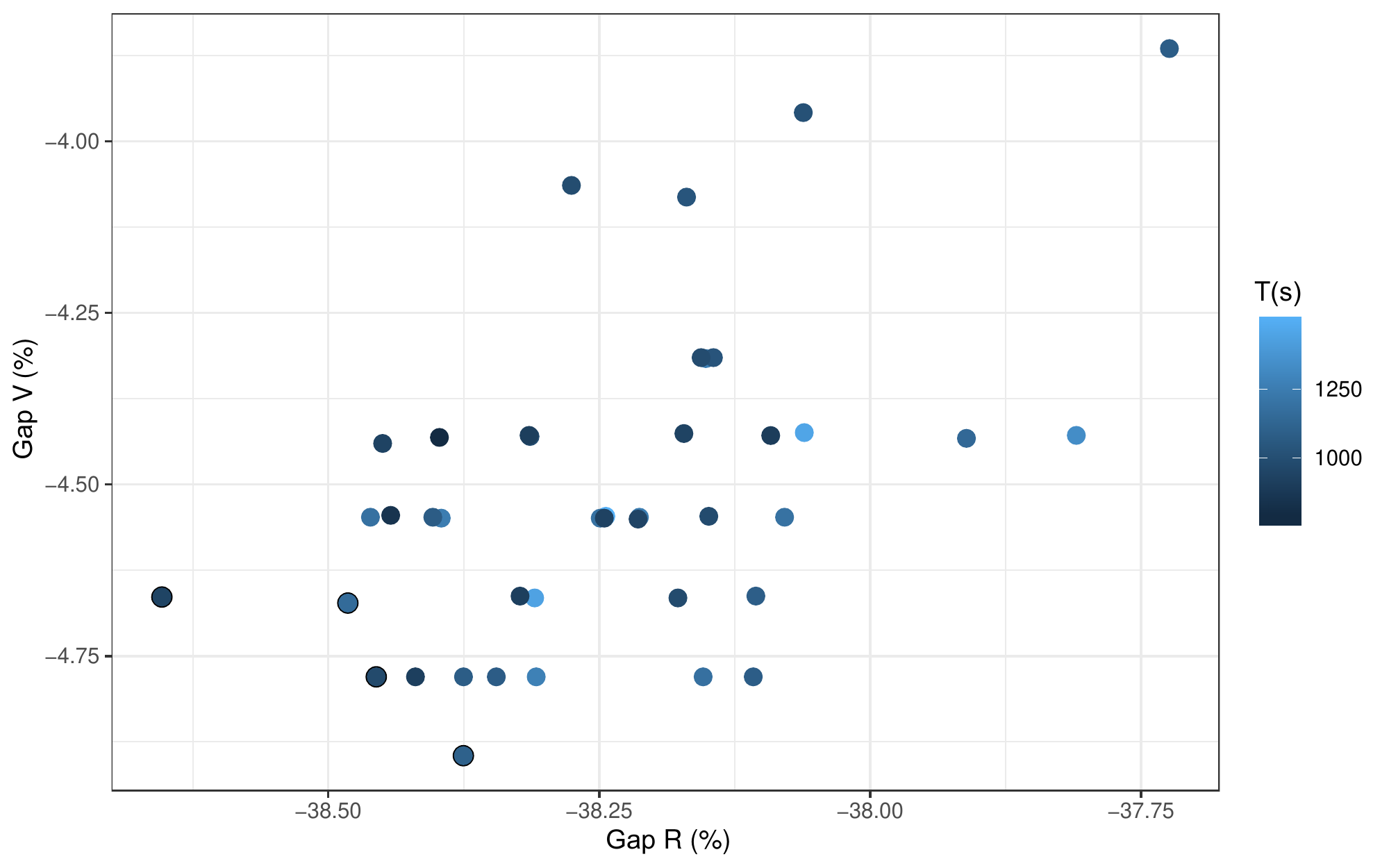}
    \caption{Sensitivity results for Diving Heuristic using 500-services instance set}
    \label{fig:diving}
\end{figure}

We report the results of the 42 configurations in Figure \ref{fig:diving}. We registered the number of vehicles and the average driving cost of the solutions of each algorithm configuration and calculated their relative gap (negative values representing improvements) from the best solutions of the ALNS. Therefore, we plot each configuration in a plane with axes representing the percentage gap in terms of driving cost and number of vehicles. The color of each dot represents the computational time of the configuration, with lighter colors indicating slower configurations.
 
Based on our experiments, \myblue{all considered configurations of the diving heuristic largely outperform the ALNS. Moreover, though relative performance differences between configurations are generally small,} two configurations stand out: the one with the best average driving cost and the one with the smallest average number of vehicles. The solution with the best driving cost uses configuration \texttt{(0, 0, 1, F)}. This configuration achieves driving-cost and fleet-size values that are $38.65\%$ and $4.66\%$ better than the ALNS solutions, respectively, for an average computational time of $920.3$ seconds. This configuration is also the simplest one, as it does not use strong diving, fixes only one route on each iteration, and does not rely on singleton routes. In contrast, the solution with the smallest average number of vehicles uses configuration \texttt{(2, 2, 1, T)}. It achieved average improvements of $38.38\%$ and $4.90\%$ in driving cost and fleet size over the ALNS solutions, for an average computational time of $1083.2$ seconds. This configuration exploits all the features discussed in Section~\ref{sec:diving}. First, it pre-generates singleton routes. Then, at each iteration, the strong diving considers two alternative routes and selects the best. This strategy is used up to a depth of two, after which the standard diving is used with one route fixed per iteration.

Since the objective value of the EVSP prioritizes fleet-size minimization on these instances, we will provide detailed results for configuration \texttt{(2, 2, 1, T)} as it is the best for this criterion. Tables~\ref{tbl:dive100} and~\ref{tbl:dive500} show the results for the instances with 100 and 500 services, respectively. The first four columns of~Table \ref{tbl:dive100} present the best ALNS solutions, and the following four columns present the optimal results found by the regular B\&P (from Table \ref{tbl:bap100}). The eight remaining columns give the results of the diving algorithm and compare them to the ALNS in terms of their upper bound (UB and associated gap), number of vehicles (V and associated gap), driving cost (D and associated gap), computational time (T) in seconds, \myblue{number of nodes fixed (\#\textsubscript{F}), and depth (\#\textsubscript{D}).}

\begin{table}[htbp]%
\centering
\setlength{\tabcolsep}{4pt}
\resizebox{\textwidth}{!}{%
\begin{tabular}{lccccccccccccccccccc}\toprule
\multirow{2}{*}{Instance} & \multicolumn{4}{c}{ALNS} && \multicolumn{4}{c}{B\&P} && \multicolumn{8}{c}{Diving Heuristic} \\ \cline{2-5}\cline{7-10}\cline{12-20}
 & \multicolumn{1}{c}{UB} & \multicolumn{1}{c}{V} & \multicolumn{1}{c}{D} & \multicolumn{1}{c}{T} && \multicolumn{1}{c}{Opt} & \multicolumn{1}{c}{V} & \multicolumn{1}{c}{D} & \multicolumn{1}{c}{T} && \multicolumn{1}{c}{UB} & \multicolumn{1}{c}{Gap} & \multicolumn{1}{c}{V} & \multicolumn{1}{c}{Gap} & \multicolumn{1}{c}{D} & \multicolumn{1}{c}{Gap} & \multicolumn{1}{c}{T} & \multicolumn{1}{c}{\#\textsubscript{F}} & \multicolumn{1}{c}{\#\textsubscript{D}} \\\midrule
\texttt{D2\_S4\_C100\_1} & 211775.0 & 21 & 1775.0 & 252 && \textbf{211741.0} & 21 & 1741.0 & 22.1 && 211788.1 & 0.01\% & 21 & 0.00\% & 1788.1 & 0.74\% & 4.4 & 16 & 13 \\
\texttt{D2\_S4\_C100\_2} & 182178.0 & 18 & 2178.0 & 890 && \textbf{181932.1} & 18 & 1932.1 & 39.2 && 192024.3 & 5.40\% & 19 & 5.56\% & 2024.3 & -7.06\% & 3.8 & 20 & 17 \\
\texttt{D2\_S4\_C100\_3} & 192230.0 & 19 & 2230.0 & 998 && \textbf{182231.7} & 18 & 2231.7 & 8.2 && 182356.3 & -5.14\% & 18 & -5.26\% & 2356.3 & 5.66\% & 3.1 & 12 & 9 \\
\texttt{D2\_S4\_C100\_4} & 212231.0 & 21 & 2231.0 & 319 && \textbf{212115.7} & 21 & 2115.7 & 96.5 && 212177.9 & -0.03\% & 21 & 0.00\% & 2177.9 & -2.38\% & 6.3 & 7 & 4 \\
\texttt{D2\_S4\_C100\_5} & 181882.0 & 18 & 1882.0 & 381 && \textbf{181685.2} & 18 & 1685.2 & 12.9 && \textbf{181685.2} & -0.11\% & 18 & 0.00\% & 1685.2 & -10.46\% & 3.8 & 7 & 4 \\
\texttt{D4\_S8\_C100\_1} & 191600.0 & 19 & 1600.0 & 281 && \textbf{191470.7} & 19 & 1470.7 & 21.6 && 191499.8 & -0.05\% & 19 & 0.00\% & 1499.8 & -6.26\% & 4.3 & 13 & 10 \\
\texttt{D4\_S8\_C100\_2} & 192097.0 & 19 & 2097.0 & 500 && \textbf{191902.5} & 19 & 1902.5 & 22.5 && 191985.5 & -0.06\% & 19 & 0.00\% & 1985.5 & -5.32\% & 4.4 & 12 & 9 \\
\texttt{D4\_S8\_C100\_3} & 191510.0 & 19 & 1510.0 & 325 && \textbf{191401.7} & 19 & 1401.7 & 153.9 && 191412.3 & -0.05\% & 19 & 0.00\% & 1412.3 & -6.47\% & 3.8 & 9 & 6 \\
\texttt{D4\_S8\_C100\_4} & 211612.0 & 21 & 1612.0 & 306 && \textbf{211468.4} & 21 & 1468.4 & 50.9 && 211597.1 & -0.01\% & 21 & 0.00\% & 1597.1 & -0.92\% & 5.5 & 15 & 12 \\
\texttt{D4\_S8\_C100\_5} & 191704.0 & 19 & 1704.0 & 311 && \textbf{191592.5} & 19 & 1592.5 & 40.1 && 191628.4 & -0.04\% & 19 & 0.00\% & 1628.4 & -4.43\% & 4.7 & 14 & 11 \\\midrule
Average & 195881.9 & 19.4 & 1881.9 & 456.3 && 194754.1 & 19.3 & 1754.1 & 46.8 && 195815.5 & -0.01\% & 19.4 & 0.03\% & 1815.5 & -3.69\% & 4.4 & 12.5 & 9.5 \\\bottomrule
\end{tabular}
}
\caption{Results of the diving heuristic for the 100-services instances}
\label{tbl:dive100}
\end{table}

As observed in Table~\ref{tbl:dive100}, the diving heuristic finds solutions of a quality comparable or better than the ALNS, with a similar number of vehicles (one more in the case of instance \texttt{D2\_S4\_C100\_2}, one less for \texttt{D2\_S4\_C100\_3}) and a 3.69\% average improvement of driving cost. It is also noteworthy that the diving heuristic is \myblue{two orders of magnitude faster} than the ALNS on those medium-scale instances, as it uses 4.4 seconds on average compared to $5 \times 456.3$ seconds for  ALNS.

Compared to the optimal B\&P solutions, the number of vehicles and driving costs of the diving heuristic are $0.56\%$ and $3.41\%$ higher on average, but the computational time is $10.6\times$ smaller ($5.6\times$ faster compared to the B\&P with sparsification). Notably, the diving heuristic found the optimal solution on instance \texttt{D2\_S4\_C100\_5}, and it found the optimal number of vehicles in all cases but \texttt{D2\_S4\_C100\_2}.

Table~\ref{tbl:dive500} now presents the results for the 500-services instance set using a similar format as Table~\ref{tbl:dive100}. As the optimal solutions are unknown for these instances, we show the best lower bound provided by the B\&P, when available.

\begin{table}[htbp]%
\centering
\setlength{\tabcolsep}{4pt}
\resizebox{\textwidth}{!}{%
\begin{tabular}{lcccccccccccccccc}\toprule
\multirow{2}{*}{Instance} & \multicolumn{4}{c}{ALNS} && \multicolumn{1}{c}{B\&P} && \multicolumn{8}{c}{Diving Heuristic} \\ \cline{2-5} \cline{7-7} \cline{9-17}
& \multicolumn{1}{c}{UB} & \multicolumn{1}{c}{V} & \multicolumn{1}{c}{D} & \multicolumn{1}{c}{T} && \multicolumn{1}{c}{LB} && \multicolumn{1}{c}{UB} & \multicolumn{1}{c}{Gap} & \multicolumn{1}{c}{V} & \multicolumn{1}{c}{Gap} & \multicolumn{1}{c}{D} & \multicolumn{1}{c}{Gap} & \multicolumn{1}{c}{T} & \multicolumn{1}{c}{\#\textsubscript{F}} & \multicolumn{1}{c}{\#\textsubscript{D}} \\\midrule
\texttt{D4\_S8\_C500\_1} & 878650.0 & 87 & 8650.0  & 1233 && 835222.8 && \textbf{835581.4} & -4.90\% & 83 & -4.60\% & 5581.4 & -35.47\% & 549.3 & 77 & 74 \\
\texttt{D4\_S8\_C500\_2} & 940142.0 & 93 & 10142.0 & 1256 && 925500.6 && \textbf{926034.3} & -1.50\% & 92 & -1.08\% & 6034.3 & -40.50\% & 334.1 & 80 & 77 \\
\texttt{D4\_S8\_C500\_3} & 859788.0 & 85 & 9788.0  & 1187 && -- && \textbf{766242.6} & -10.88\% & 76 & -10.59\% & 6242.6 & -36.22\% & 2472.5 & 67 & 64 \\
\texttt{D4\_S8\_C500\_4} & 870033.0 & 86 & 10033.0 & 1187 && 825646.6 && \textbf{826230.3} & -5.03\% & 82 & -4.65\% & 6230.3 & -37.90\% & 1934.3 & 77 & 64 \\
\texttt{D4\_S8\_C500\_5} & 880386.0 & 87 & 10386.0 & 1198 && 786318.8 && \textbf{787294.8} & -10.57\% & 78 & -10.34\% & 7294.8 & -29.76\% & 2778.7 & 62 & 59 \\
\texttt{D8\_S16\_C500\_1} & 869530.0 & 86 & 9530.0 & 1255 && 814907.3 && \textbf{815278.9} & -6.24\% & 81 & -5.81\% & 5278.9 & -44.61\% & 462.9 & 70 & 67 \\
\texttt{D8\_S16\_C500\_2} & 869282.0 & 86 & 9282.0 & 1259 && 824939.6 && \textbf{825389.7} & -5.05\% & 82 & -4.65\% & 5389.7 & -41.93\% & 669.0 & 68 & 65 \\
\texttt{D8\_S16\_C500\_3} & 877456.0 & 87 & 7456.0 & 1178 && 874670.3 && \textbf{875060.7} & -0.27\% & 87 & 0.00\% & 5060.7 & -32.13\% & 755.9 & 80 & 77 \\
\texttt{D8\_S16\_C500\_4} & 828538.0 & 82 & 8538.0 & 1189 && 784539.5 && \textbf{784832.4} & -5.28\% & 78 & -4.88\% & 4832.4 & -43.40\% & 263.6 & 67 & 64 \\
\texttt{D8\_S16\_C500\_5} & 858816.0 & 85 & 8816.0 & 1183 && -- && \textbf{835129.1} & -2.76\% & 83 & -2.35\% & 5129.1 & -41.82\% & 611.6 & 71 & 68 \\\midrule
Average & 873262.1 & 86.4 & 9262.1 & 1212.5 && -- && 827707.4 & -5.25\% & 82.2 & -4.90\% & 5707.4 & -38.38\% & 1083.2 & 71.9 & 68.9 \\\bottomrule
\end{tabular}
}
\caption{Diving heuristic results for 500-services instance set}
\label{tbl:dive500}
\end{table}

As visible on these experiments, the diving heuristic vastly surpasses the ALNS in terms of solution quality on these large-scale instances, thereby providing a meaningful solution alternative for these challenging cases. We measure an average driving-cost improvement of $38.38\%$,  suggesting that the ALNS was primarily calibrated for fleet-size minimization. In terms of fleet size, the diving heuristic produced solutions with fewer vehicles than the ALNS for 9 out of the 10 instances, corresponding to a 4.90\% reduction in the number of vehicles on average.

Noticeably, a comparison with the available lower bounds produced by the regular B\&P algorithm (without sparsification) shows that the diving heuristic solutions are no further than $0.06\%$ from the optimal solution cost. With these values, the number of vehicles found by the diving heuristics is very likely to be optimum in all cases. Finally, the heuristic solution took $1083.2$ seconds on average, whereas the B\&P was interrupted after 21,600 seconds of computational effort without an optimal solution. It is important to mention that throughout the diving heuristic experiments, there was no case of infeasibility using the proposed parameters. \\

\myred{Finally, to analyze more finely the impact of the number of requests, we tested the diving heuristic on the same instances with 200, 300, and 400 services as in Section~\ref{sec:BPresults}.
On these instances, we measured average optimality gaps between the best lower bound found by the B\&P on the non-sparsified instances and the primal solution of the diving heuristic of $0.015\%$, $0.011\%$, and $0.019\%$, respectively. The diving algorithm used 31.5, 100.8, and 368.0 seconds of computational time on these instances, respectively. Therefore, it appears to scale very well in terms of solution quality, with a computational time growth that is slightly worse than cubic in the number of requests.}

\section{Conclusions}
\label{sec:conclusion}

In this paper, we have introduced efficient branch-and-price and diving algorithms for the EVSP. \myblue{Our methods exploit an efficient characterization of non-dominated charging arcs, along with backward bounds within a bi-directional pricing algorithm. We also developed sparsification techniques and diving strategies to quickly find high-quality primal solutions for large instances.} Thanks to this methodology, we could find optimal solutions for all 100-services instances with the B\&P algorithm, and solutions no further than $0.06\%$ from the optimum with the diving heuristic for the challenging 500-services instances. Compared to previous solution approaches, the proposed diving heuristic finds much higher quality solutions in a shorter time.

The research perspectives connected to this work are numerous. First of all, research can be pursued to increase the scalability of the proposed algorithms. An option along this line is to solve the resource-constrained shortest path problems heuristically, possibly by defining a stronger heuristic dominance. Another promising research perspective concerns the integration of machine learning components for critical steps of the method, e.g., for the sparsification step \citep[e.g., as in][]{Joshi2019} or the branch choices in the B\&P and diving algorithms \citep[see, e.g.,][]{Gasse2019,BengioSurvey2021}. Finally, we recommend pursuing the study of pricing algorithms based on monoid structures as done in \citep{Parmentier2018}, considering a broader range of applications to routing and scheduling problems with non-linear extension functions and evaluating their performance in dynamic and stochastic settings.

\section*{Acknowledgments}

The authors would like to thank Stefan R{\o}pke for his assistance with the instance files and experimental conventions, as well as Victor Abu-Marrul for help on Figure~\ref{fig:diving}. This research was partially supported by CAPES with Finance Code 001, by CNPq under grants 308528/2018-2 and 315361/2020-4, and by FAPERJ under grants E-26/202.790/2019 and E-26/010.002232/2019 in Brazil. This financial support is gratefully acknowledged.

\bibliographystyle{ormsv080-noURLDOI}
\bibliography{references}

\begin{thebibliography}{47}
\expandafter\ifx\csname natexlab\endcsname\relax\def\natexlab#1{#1}\fi
\expandafter\ifx\csname url\endcsname\relax
  \def\url#1{{\tt #1}}\fi
\expandafter\ifx\csname urlprefix\endcsname\relax\def\urlprefix{URL }\fi
\expandafter\ifx\csname urlstyle\endcsname\relax
  \expandafter\ifx\csname doi\endcsname\relax
  \def\doi#1{doi:\discretionary{}{}{}#1}\fi \else
  \expandafter\ifx\csname doi\endcsname\relax
  \def\doi{doi:\discretionary{}{}{}\begingroup \urlstyle{rm}\Url}\fi \fi

\bibitem[{Adler and Mirchandani(2017)}]{Adler2017}
Adler, J.D., P.B. Mirchandani. 2017.
\newblock {The vehicle scheduling problem for fleets with alternative-fuel
  vehicles}.
\newblock {\it Transportation Science\/} {\bf 51}(2) 441--456.

\bibitem[{Barnhart et~al.(1998)Barnhart, Johnson, Nemhauser, Savelsbergh, and
  Vance}]{Barnhart1998}
Barnhart, C., E.L. Johnson, G.L. Nemhauser, M.W.P. Savelsbergh, P.H. Vance.
  1998.
\newblock {Branch-and-price: Column generation for solving huge integer
  programs}.
\newblock {\it Operations Research\/} {\bf 46}(3) 316--329.

\bibitem[{Bengio et~al.(2021)Bengio, Lodi, and Prouvost}]{BengioSurvey2021}
Bengio, Y., A.~Lodi, A.~Prouvost. 2021.
\newblock Machine learning for combinatorial optimization: A methodological
  tour d’horizon.
\newblock {\it European Journal of Operational Research\/} {\bf 290}(2)
  405--421.

\bibitem[{BloombergNEF(2021)}]{BloombergNEF2021}
BloombergNEF. 2021.
\newblock {Electric Vehicle Outlook}.
\newblock Tech. rep.

\bibitem[{Bodin et~al.(1983)Bodin, Golden, Assad, and Ball}]{Bodin1983}
Bodin, L., B.~Golden, A.~Assad, M.~Ball. 1983.
\newblock {Routing and scheduling of vehicles and crews: The state of the art.}
\newblock {\it Computers {\&} Operations Research\/} {\bf 10}(2) 63--211.

\bibitem[{Bol{\'{i}}var et~al.(2014)Bol{\'{i}}var, Lozano, and
  Medaglia}]{Bolivar2014}
Bol{\'{i}}var, M.A., L.~Lozano, A.L. Medaglia. 2014.
\newblock {Acceleration strategies for the weight constrained shortest path
  problem with replenishment}.
\newblock {\it Optimization Letters\/} {\bf 8}(8) 2155--2172.

\bibitem[{Brandst{\"{a}}tter et~al.(2016)Brandst{\"{a}}tter, Gambella, Leitner,
  Malaguti, Masini, Puchinger, Ruthmair, and Vigo}]{Brandstatter2016}
Brandst{\"{a}}tter, G., C.~Gambella, M.~Leitner, E.~Malaguti, Fi. Masini,
  J.~Puchinger, M.~Ruthmair, D.~Vigo. 2016.
\newblock {Overview of optimization problems in electric car-sharing system
  design and management}.
\newblock H.~Dawid, K.F. Doerner, G.~Feichtinger, P.M. Kort, A.~Seidl, eds.,
  {\it Dynamic Perspectives on Managerial Decision Making\/}. Springer
  International Publishing, 441--471.

\bibitem[{Breunig et~al.(2019)Breunig, Baldacci, Hartl, and
  Vidal}]{Breunig2018}
Breunig, U., R.~Baldacci, R.~Hartl, T.~Vidal. 2019.
\newblock {The electric two-echelon vehicle routing problem}.
\newblock {\it Computers {\&} Operations Research\/} {\bf 103} 198--210.

\bibitem[{Bunte and Kliewer(2009)}]{Bunte2009}
Bunte, S., N~Kliewer. 2009.
\newblock {An overview on vehicle scheduling models}.
\newblock {\it Public Transport\/} {\bf 1}(4) 299--317.

\bibitem[{Cabrera et~al.(2020)Cabrera, Medaglia, Lozano, and
  Duque}]{Cabrera2020}
Cabrera, N., A.L. Medaglia, L.~Lozano, D.~Duque. 2020.
\newblock {An exact bidirectional pulse algorithm for the constrained shortest
  path}.
\newblock {\it Networks\/} {\bf 76}(2) 128--146.

\bibitem[{{De Jong}(2018)}]{DeJong2018}
{De Jong}, S. 2018.
\newblock {London's iconic black cabs go electric}.

\bibitem[{Desaulniers et~al.(2016)Desaulniers, Errico, Irnich, and
  Schneider}]{Desaulniers2016}
Desaulniers, G., F.~Errico, S.~Irnich, M.~Schneider. 2016.
\newblock {Exact algorithms for electric vehicle-routing problems with time
  windows}.
\newblock {\it Operations Research\/} {\bf 64}(6) 1388--1405.

\bibitem[{Desrosiers et~al.(1995)Desrosiers, Dumas, Solomon, and
  Soumis}]{Desrosiers1995}
Desrosiers, J., Y.~Dumas, M.M. Solomon, F.~Soumis. 1995.
\newblock {Time constrained routing and scheduling}.
\newblock M.~Ball, T.~L. Magnanti, C.L. Monma, G.L. Nemhauser, eds., {\it
  Network Routing\/}. North-Holland Amsterdam, 35--139.

\bibitem[{Dumitrescu and Boland(2003)}]{Dumitrescu2003}
Dumitrescu, I., N.~Boland. 2003.
\newblock {Improved preprocessing, labeling and scaling algorithms for the
  weight-constrained shortest path problem}.
\newblock {\it Networks\/} {\bf 42}(3) 135--153.

\bibitem[{Dunne(2017)}]{Dunne2017}
Dunne, M.J. 2017.
\newblock {China deploys aggressive mandates to take lead in electric
  vehicles}.

\bibitem[{Felipe et~al.(2014)Felipe, Ortu{\~{n}}o, Righini, and
  Tirado}]{Felipe2014}
Felipe, {\'{A}}., M.T. Ortu{\~{n}}o, G.~Righini, G.~Tirado. 2014.
\newblock {A heuristic approach for the green vehicle routing problem with
  multiple technologies and partial recharges}.
\newblock {\it Transportation Research Part E: Logistics and Transportation
  Review\/} {\bf 71} 111--128.

\bibitem[{Gasse et~al.(2019)Gasse, Chetelat, Ferroni, Charlin, and
  Lodi}]{Gasse2019}
Gasse, M., D.~Chetelat, N.~Ferroni, L.~Charlin, A.~Lodi. 2019.
\newblock Exact combinatorial optimization with graph convolutional neural
  networks.
\newblock H.~Wallach, H.~Larochelle, A.~Beygelzimer, F.~d'Alch\'{e} Buc,
  E.~Fox, R.~Garnett, eds., {\it Advances in Neural Information Processing
  Systems\/}, vol.~32. Curran Associates, Inc.

\bibitem[{Haghani and Banihashemi(2002)}]{Haghani2002}
Haghani, A., M.~Banihashemi. 2002.
\newblock {Heuristic approaches for solving large-scale bus transit vehicle
  scheduling problem with route time constraints}.
\newblock {\it Transportation Research Part A\/} {\bf 36} 309--333.

\bibitem[{Hiermann et~al.(2019)Hiermann, Hartl, Puchinger, and
  Vidal}]{Hiermann2019}
Hiermann, G., R.~Hartl, J.~Puchinger, T.~Vidal. 2019.
\newblock {Routing a mix of conventional, plug-in hybrid, and electric
  vehicles}.
\newblock {\it European Journal of Operational Research\/} {\bf 272}(1)
  235--248.

\bibitem[{Hiermann et~al.(2016)Hiermann, Puchinger, Ropke, and
  Hartl}]{Hiermann2016}
Hiermann, G., J.~Puchinger, S.~Ropke, R.F. Hartl. 2016.
\newblock {The electric fleet size and mix vehicle routing problem with time
  windows and recharging stations}.
\newblock {\it European Journal of Operational Research\/} {\bf 252}(3)
  995--1018.

\bibitem[{Hu et~al.(2018)Hu, Dong, Lin, and Yang}]{Hu2018}
Hu, L., J.~Dong, Z.~Lin, J.~Yang. 2018.
\newblock {Analyzing battery electric vehicle feasibility from taxi travel
  patterns: The case study of New York City, USA}.
\newblock {\it Transportation Research Part C: Emerging Technologies\/} {\bf
  87} 91--104.

\bibitem[{Irnich and Desaulniers(2005)}]{Irnich2005}
Irnich, S., G.~Desaulniers. 2005.
\newblock {Shortest path problems with resource constraints}.
\newblock G.~Desaulniers, J.~Desrosiers, M.M. Solomon, eds., {\it Column
  Generation\/}. Springer, New York, 33--65.

\bibitem[{Joshi et~al.(2019)Joshi, Laurent, and Bresson}]{Joshi2019}
Joshi, C.K., T.~Laurent, X.~Bresson. 2019.
\newblock An efficient graph convolutional network technique for the travelling
  salesman problem.
\newblock {\it arXiv preprint arXiv:1906.01227\/} .

\bibitem[{Kergosien et~al.(2021)Kergosien, Giret, Neron, and
  Sauvanet}]{kergosien2021efficient}
Kergosien, Yannick, Antoine Giret, Emmanuel Neron, Ga{\"e}l Sauvanet. 2021.
\newblock An efficient label-correcting algorithm for the multiobjective
  shortest path problem.
\newblock {\it INFORMS Journal on Computing\/} .

\bibitem[{Keskin and {\c{C}}atay(2018)}]{Keskin2018}
Keskin, M., B.~{\c{C}}atay. 2018.
\newblock {A matheuristic method for the electric vehicle routing problem with
  time windows and fast chargers}.
\newblock {\it Computers {\&} Operations Research\/} {\bf 100} 172--188.

\bibitem[{Li(2014)}]{Li2014}
Li, J.-Q. 2014.
\newblock {Transit bus scheduling with limited energy}.
\newblock {\it Transportation Science\/} {\bf 48}(4) 521--539.

\bibitem[{Li et~al.(2016)Li, Zhan, de~Jong, and Lukszo}]{Li2016c}
Li, Y., C.~Zhan, M.~de~Jong, Z.~Lukszo. 2016.
\newblock {Business innovation and government regulation for the promotion of
  electric vehicle use: lessons from Shenzhen, China}.
\newblock {\it Journal of Cleaner Production\/} {\bf 134} 371--383.

\bibitem[{{London Authorities}(2018)}]{London2018}
{London Authorities}. 2018.
\newblock {Mayor marks key milestone of 100 rapid charging points across
  London}.

\bibitem[{Mahmoud et~al.(2016)Mahmoud, Garnett, Ferguson, and
  Kanaroglou}]{Mahmoud2016}
Mahmoud, M., R.~Garnett, M.~Ferguson, P.~Kanaroglou. 2016.
\newblock {Electric buses: A review of alternative powertrains}.
\newblock {\it Renewable and Sustainable Energy Reviews\/} {\bf 62} 673--684.

\bibitem[{Martinelli et~al.(2011)Martinelli, Pecin, Poggi, and
  Longo}]{martinelli2011branch}
Martinelli, R., D.~Pecin, M.~Poggi, H.~Longo. 2011.
\newblock A branch-cut-and-price algorithm for the capacitated arc routing
  problem.
\newblock {\it International Symposium on Experimental Algorithms\/}. Springer,
  315--326.

\bibitem[{Montoya et~al.(2017)Montoya, Gu{\'{e}}ret, Mendoza, and
  Villegas}]{Montoya2017}
Montoya, A., C.~Gu{\'{e}}ret, J.E. Mendoza, J.G. Villegas. 2017.
\newblock {The electric vehicle routing problem with nonlinear charging
  function}.
\newblock {\it Transportation Research Part B: Methodological\/} {\bf 103}
  87--110.

\bibitem[{Parmentier(2019)}]{Parmentier2018}
Parmentier, A. 2019.
\newblock {Algorithms for non-linear and stochastic resource constrained
  shortest paths}.
\newblock {\it Mathematical Methods of Operations Research\/} {\bf 89}
  281--317.

\bibitem[{Pelletier et~al.(2016)Pelletier, Jabali, and Laporte}]{Pelletier2016}
Pelletier, S., O.~Jabali, G.~Laporte. 2016.
\newblock {50th anniversary invited article---Goods distribution with electric
  vehicles: Review and research perspectives}.
\newblock {\it Transportation Science\/} {\bf 50}(1) 3--22.

\bibitem[{Pelletier et~al.(2017)Pelletier, Jabali, Laporte, and
  Veneroni}]{Pelletier2017}
Pelletier, S., O.~Jabali, G.~Laporte, M.~Veneroni. 2017.
\newblock {Battery degradation and behaviour for electric vehicles: Review and
  numerical analyses of several models}.
\newblock {\it Transportation Research Part B: Methodological\/} {\bf 103}
  158--187.

\bibitem[{Pessoa et~al.(2010)Pessoa, Uchoa, De~Arag{\~a}o, and
  Rodrigues}]{pessoa2010exact}
Pessoa, Artur, Eduardo Uchoa, Marcus~Poggi De~Arag{\~a}o, Rosiane Rodrigues.
  2010.
\newblock Exact algorithm over an arc-time-indexed formulation for parallel
  machine scheduling problems.
\newblock {\it Mathematical Programming Computation\/} {\bf 2}(3) 259--290.

\bibitem[{Ribeiro and Soumis(1994)}]{Ribeiro1994}
Ribeiro, C.C., F.~Soumis. 1994.
\newblock {A column generation approach to the multiple-depot vehicle
  scheduling problem}.
\newblock {\it Operations Research\/} {\bf 42}(1) 41--52.

\bibitem[{Sadykov et~al.(2019)Sadykov, Vanderbeck, Pessoa, Tahiri, and
  Uchoa}]{sadykov2019primal}
Sadykov, R., F.~Vanderbeck, A.~Pessoa, I.~Tahiri, E.~Uchoa. 2019.
\newblock Primal heuristics for branch and price: The assets of diving methods.
\newblock {\it INFORMS Journal on Computing\/} {\bf 31}(2) 251--267.

\bibitem[{Schiffer et~al.(2019)Schiffer, Schneider, Walther, and
  Laporte}]{Schiffer2018}
Schiffer, M., M.~Schneider, G.~Walther, G.~Laporte. 2019.
\newblock {Vehicle routing and location-routing with intermediate stops: A
  review}.
\newblock {\it Transportation Science\/} {\bf 53}(2) 319--343.

\bibitem[{Schiffer and Walther(2017)}]{Schiffer2017}
Schiffer, M., G.~Walther. 2017.
\newblock {The electric location routing problem with time windows and partial
  recharging}.
\newblock {\it European Journal of Operational Research\/} {\bf 260}(3)
  995--1013.

\bibitem[{Scorrano et~al.(2020)Scorrano, Danielis, and
  Giansoldati}]{Scorrano2020}
Scorrano, M., R.~Danielis, M.~Giansoldati. 2020.
\newblock {Mandating the use of the electric taxis: The case of Florence}.
\newblock {\it Transportation Research Part A: Policy and Practice\/} {\bf 132}
  402--414.

\bibitem[{Smith et~al.(2012)Smith, Boland, and Waterer}]{Smith2012}
Smith, O.J., N.~Boland, H.~Waterer. 2012.
\newblock {Solving shortest path problems with a weight constraint and
  replenishment arcs}.
\newblock {\it Computers {\&} Operations Research\/} {\bf 39}(5) 964--984.

\bibitem[{Thomas et~al.(2019)Thomas, Calogiuri, and Hewitt}]{Thomas2019}
Thomas, B.W., T.~Calogiuri, M.~Hewitt. 2019.
\newblock {An exact bidirectional A* approach for solving resource-constrained
  shortest path problems}.
\newblock {\it Networks\/} {\bf 73}(2) 187--205.

\bibitem[{Tilk et~al.(2017)Tilk, Rothenb{\"{a}}cher, Gschwind, and
  Irnich}]{Tilk2017}
Tilk, C., A.-K. Rothenb{\"{a}}cher, T.~Gschwind, S.~Irnich. 2017.
\newblock {Asymmetry matters: Dynamic half-way points in bidirectional labeling
  for solving shortest path problems with resource constraints faster}.
\newblock {\it European Journal of Operational Research\/} {\bf 261}(2)
  530--539.

\bibitem[{{van Kooten Niekerk} et~al.(2017){van Kooten Niekerk}, van~den Akker,
  and Hoogeveen}]{VanKootenNiekerk2017}
{van Kooten Niekerk}, M.E., J.M. van~den Akker, J.A. Hoogeveen. 2017.
\newblock {Scheduling electric vehicles}.
\newblock {\it Public Transport\/} {\bf 9}(1-2) 155--176.

\bibitem[{Vidal et~al.(2020)Vidal, Laporte, and Matl}]{Vidal2020}
Vidal, T., G.~Laporte, P.~Matl. 2020.
\newblock {A concise guide to existing and emerging vehicle routing problem
  variants}.
\newblock {\it European Journal of Operational Research\/} {\bf 286} 401--416.

\bibitem[{Wang et~al.(2017)Wang, Huang, Xu, and Barclay}]{Wang2017c}
Wang, Y., Y.~Huang, J.~Xu, N.~Barclay. 2017.
\newblock {Optimal recharging scheduling for urban electric buses: A case study
  in Davis}.
\newblock {\it Transportation Research Part E: Logistics and Transportation
  Review\/} {\bf 100} 115--132.

\bibitem[{Wen et~al.(2016)Wen, Linde, Ropke, Mirchandani, and Larsen}]{Wen2016}
Wen, M., E.~Linde, S.~Ropke, P.~Mirchandani, A.~Larsen. 2016.
\newblock {An adaptive large neighborhood search heuristic for the Electric
  Vehicle Scheduling Problem}.
\newblock {\it Computers {\&} Operations Research\/} {\bf 76} 73--83.

\end{thebibliography}

\renewcommand\thesubsection{\Alph{subsection}}
\section*{\myblue{Appendix -- Proofs}}

We start by introducing some lemmas on~$\tch$ and $\preceq$ that will permit us to simplify our proofs of the results of Section~\ref{sec:StationSequenceScheduling}.

\subsection{Lemma on $\tch$}

\begin{lemma}
\label{lem:tchComposition}
\begin{enumerate}
        \item Given $\ell_1$, $\ell_2$, and $\ell_3$ in $[0,\mc]$, we have $\tch(\ell_1,\ell_3) = \tch(\ell_1,\ell_2) + \tch(\ell_2,\ell_3)$.
        \item Given $e \leq \ell \leq \ell'$, we have $\tch(\ell-e,\ell'-e) \leq \tch(\ell,\ell')$.
\end{enumerate}
\end{lemma}
\proof{Proof.}
Let us prove the first point for $\ell_1 \leq \ell_2 \leq \ell _3$.
Since $\tau \mapsto \varphi(\ell_1,\tau)$ is non-decreasing, we have $\tch(\ell_1,\ell_3) = \tch(\ell_1,\ell_2)$. Next we have:
\begin{align*}
    \tch(\ell_1,\ell_3) 
    &= \min\big\{\tau \colon \varphi(\ell_1,\tau) = \ell_3\} \\
    &= \tch(\ell_1,\ell_2) + \min \big\{\tau' \colon \varphi(\ell_1, \tch(\ell_1,\ell_2) + \tau') = \ell_3\} & & \text{replacing $\tau$ by $\tau' = \tau - \tch(\ell_1,\ell_2)$} \\
    &= \tch(\ell_1,\ell_2) + \underbrace{\min \big\{\tau' \colon \varphi\big(\underbrace{\varphi(\ell_1, \tch(\ell_1,\ell_2)}_{\ell_2}, \tau')\big) = \ell_3\}}_{\tch(\ell_2,\ell_3)} & &\text{using } \varphi(\ell,\tau + \tau') = \varphi\big(\varphi(\ell,\tau), \tau'\big).
\end{align*}
The other cases follow immediately for $\varphi(\ell,\ell') = -\varphi(\ell',\ell)$.
Indeed, suppose for instance  that $\ell_2\leq \ell_1 \leq \ell_3$.
Applying the result already proven, we have $\tch(\ell_2,\ell_3) = \tch(\ell_2,\ell_1) + \tch(\ell_1,\ell_3) = -\tch(\ell_1,\ell_2) + \tch(\ell_1,\ell_3)$, which gives the result.

The second point follows from the fact that $\tau \mapsto \varphi(\ell,\tau)$ is non-increasing.
Indeed, $\varphi(\ell-e,\tch(\ell - \ell')) \geq \ell-e + \varphi(\ell,\tch(\ell - \ell')) \geq \ell' - e$.
\halmos\endproof

\subsection{Lemmas on $\preceq$}
\label{sub:preorder}

Recall that $S = \big([0,M] \cup \{-\infty\}\big) \times \bbR$,
and $\preceq$ is the binary relation defined by:
$$ 
\begin{cases}
    (-\infty,t) \preceq (\ell',t'), & \text{for any }(\ell',t') \in S \quad \text{ and $t$ in $\bbR$}, \\ 
    (\ell,t) \preceq (-\infty,t')  \quad \text{if and only if} \quad \ell = -\infty, & \text{for any $t,t'$ in $\bbR$}, \\
    (\ell,t) \preceq (\ell',t') \quad \text{if and only if} \quad 
    t + \tch(\ell,\ell') \geq t', & \text{for $\ell,\ell' \in [0,\mc]$ and $t,t'$ in $\bbR$}.
\end{cases}
$$
As an immediate consequence of the definitions of $\tch$ and $\preceq$, we can replace the third case by:
    $$(\ell,t) \preceq (\ell',t') \text{if and only if} \quad 
    t \geq t' + \tch(\ell',\ell), \text{ for $\ell,\ell' \in [0,\mc]$ and $t,t'$ in $\bbR$}.$$

\begin{lemma}\label{lem:preceqIsAPreorder}
    $\preceq$ is a preorder on $S$.
\end{lemma}
\proof{Proof.}
$\varphi(\ell,\tau) = \ell$ implies that $\preceq$ is reflexive. Let us now prove that it is transitive. Let $(\ell_1,t_1)$, $(\ell_2,t_2)$, and $(\ell_3,t_3)$ be such that $(\ell_1,t_1) \preceq (\ell_2,t_2)$  and $(\ell_2,t_2)\preceq(\ell_3,t_3)$. 
If one $\ell_i$ is equal to $-\infty$, the result is immediate. Suppose now that $\ell_1,\ell_2,\ell_3 \in [0,\mc]$. 
Using Lemma~\ref{lem:tchComposition}.1, we have:
$$t_1 + \tch(\ell_1,\ell_3) = t_1 + \tch(\ell_1,\ell_2) + \tch(\ell_2,\ell_3) \geq  t_2 + \tch(\ell_2,\ell_3) \geq t_3, $$
which gives $(\ell_1,t_1) \preceq (\ell_3,t_3)$ and the result.
\halmos\endproof

Let $\equiv$ be the equivalence relation defined by $s \equiv s'$ if and only if $s\preceq s'$ and $s'\preceq s$. Remark that $\preceq$ is a total order on the quotient space $S/\equiv$.
Given $\delta$ and $e$ in $\bbR^+$ and $(t,\ell)$ in $S$, let us define:
$$ \varsigma_{\delta,e}(\ell,t) = \begin{cases}
(\ell -e,t+\delta) & \text{if }\ell \geq e \\
-\infty & \text{otherwise}.
\end{cases}  $$

\begin{lemma}
\label{lem:partialOrderAndVarphi}

\begin{enumerate}
    \item $(\ell,t) \preceq (\ell',t') \Leftrightarrow (\ell,t + \delta) \preceq (\ell',t' + \delta) $.
    \item If $e \leq \ell \leq \ell'$, then $(\ell',t')\preceq(\ell,t)$ implies $\varsigma_{\delta,e}(\ell',t')\preceq\varsigma_{\delta,e}(\ell,t)$.
    \item Given $(\ell,t) \in S$ and $\tau \geq 0$, we have $\big(\varphi(\ell,\tau),t+\tau\big) \preceq \big(\ell,t\big)$.
    \item Given $\ell,t,\tau, \overline \ell, \overline t, e$, and $\delta$, let 
    $$\overline{\tau} = \begin{cases}
    0 & \text{if } \overline \ell \geq e, \\
    \tch(\overline{\ell}, e) & \text{otherwise.}
    \end{cases} $$
    Then $(\ell,t) \preceq (\overline \ell, \overline t)$ implies $\varsigma_{e,\delta}\big(\varphi(\ell,\tau),t+\tau\big) \preceq \varsigma_{e,\delta}\big(\varphi(\overline\ell,\overline\tau),\overline t+ \overline \tau\big)$.
    \item $(\ell,t)\preceq(\ell',t')$ and $\tilde t \geq \max(t,t')$ implies that $(\varphi(\ell,\tilde t - t), \tilde t) \preceq (\varphi(\ell',\tilde t - t'), \tilde t)$.
\end{enumerate}
\end{lemma}
\proof{Proof.}
The first point is immediate given the definition of $\preceq$.

Let us now prove the second point. Since $e \leq \ell \leq \ell'$, Lemma~\ref{lem:tchComposition}.2 gives
$$ t + \delta +\tch(\ell-e,\ell'-e) \leq t + \delta + \tch(\ell,\ell') \leq t' + \delta,$$
which gives the result.

For the third point, remark that
$\tau \geq \min\big\{\tau' \colon \varphi(\ell,\tau') = \varphi(\ell,\tau)\big\} $,
which gives $ t + \tch(\ell,\varphi(\ell,\tau)) \leq \varphi(\ell,\tau)$ and the result. 

For the fourth point, if $\varphi(\ell,\tau) < e$, we have $\varsigma_{e,\delta}\big(\varphi(\ell,\tau),t+\tau\big) = -\infty$ and the result is immediate.
Suppose now that $\varphi(\ell,\tau) \geq e$. 
We first remark that the definition of $\overline{\tau}$ implies $(\overline \ell,\overline t) \equiv \big(\varphi(\overline\ell,\overline\tau),\overline t+ \overline \tau\big)$. Indeed, if $\overline\ell \geq e$, $(\overline \ell,\overline t) = \big(\varphi(\overline\ell,\overline\tau),\overline t+ \overline \tau\big)$, and if $\overline\ell < e$, it follows immediately from the definition of $\overline \tau$ and $\preceq$.
Combining this result with the third point of the lemma, we get  
$$ \big(\underbrace{\varphi(\ell,\tau)}_{\geq e},t+\tau\big) \preceq (\ell,t) \preceq (\overline \ell,\overline t) \equiv \big(\varphi(\overline\ell,\overline\tau),\overline t+ \overline \tau\big). $$
Applying the second point with $\big(\underbrace{\varphi(\ell,\tau)}_{\geq e},t+\tau\big)$ and $\big(\varphi(\overline\ell,\overline\tau),\overline t+ \overline \tau\big)$ then immediately gives the result.

Let $\ell$, $t$, $\ell'$, $t'$ and $\tilde t$ be as in the hypothesis of the fifth point.
We only have to prove that $\varphi(\ell,\tilde t - t) \leq \varphi(\ell', \tilde t - t') $.
Let $\overline{\ell} = \max(\ell,\ell')$, $ \overline{t} = t + \tch(\ell, \overline \ell )$, and $ \overline{t'} = t' + \tch(\ell', \overline \ell)$.
Inequality $(\ell,t) \preceq (\ell',t')$ then  gives
$ \overline{t} = t + \tch(\ell, \overline \ell) = t + \tch(\ell,\ell') + \tch(\ell',\overline \ell) \geq t' + \tch(\ell',\overline \ell) = \overline{t'}$.
Hence:
\begin{align*}
    \varphi(\ell',\tilde t - t') 
    &= \varphi(\ell',\tilde t - \overline t + \overline t - \overline{t'} + \overline{t'} - t') \\
    &= \varphi\big(\underbrace{\varphi\big(\underbrace{\varphi\big(\ell', \overline{t'} - t' \big)}_{\overline  \ell }, \underbrace{\overline t - \overline{t'} }_{\geq 0} \big)}_{\geq \overline{\ell}},\tilde t - \overline t\big) \\
    &\geq \varphi\big(\overline{\ell}, \tilde t - \overline t\big) 
    =\varphi\big(\varphi(\ell, \overline{t} - t), \tilde t - \overline t\big) 
    = \varphi(\ell, \tilde t - t),
\end{align*}
which gives the result.
\halmos
\endproof

\subsection{Proof of Theorem~\ref{prop:optimalSchedulingOfStationSequence}}

Let $a$ be the station sequence $u=s_0,s_1,\ldots,s_k,s_{k+1}=v$ between $u$ and $v$.
Let $\btr[a]$ be a scheduling of route $a$ satisfying~\eqref{eq:stationSequenceTimes}, and $\lin$ in $[0,M]\cup\{-\infty\}$ be the charge level at the end of $u$.
Let us prove by induction on $i$ from $1$ to $k$ that:
\begin{equation}
    \label{eq:optimalSchedulingProofInductionHypothesis}
    \big(\levb[s_i][a](\lin,\btr[a]),\tb[s_i][a]\big)  \preceq \big( \levb[s_i][a,*](\lin),\tb[s_i][a,*]\big).
\end{equation}
For $i=1$, since $\btr[a]$ satisfies~\eqref{eq:stationSequenceTimes}, and by definition of~$\btr[a,*](\lin)$, we have $\tb[s_1][a,*](\lin) = \tb[s_1][a] = \te[u] + \ti[u,s_1]$, and
 $\levb[s_1][a,*](\lin) = \levb[s_1][a](\lin,\btr[a]) =  \pos(\lin - \fu[u,s_1])$.
Since the preorder $\preceq$ is reflexive, we get the induction hypothesis~\eqref{eq:optimalSchedulingProofInductionHypothesis} for $i=1$.

Suppose that~\eqref{eq:optimalSchedulingProofInductionHypothesis} is true up to $i$.
Defining
$\overline{\tau} = \begin{cases}
    0 & \text{if } \levb[s_i][a,*](\lin) \geq e_{s_i,s_{i+1}} \\
    \tch(\levb[s_i][a,*](\lin), e_{s_i,s_{i+1}}) & \text{otherwise.}
    \end{cases} $
    and  $\tau = \te[s_i][a] - \tb[s_i][a]$, we get:
\begin{align*}
    \big(\levb[s_{i+1}][a](\lin,\btr[a]),\tb[s_{i+1}][a]\big) &= \varsigma_{e_{s_i,s_{i+1}},\delta_{s_i,s_{i+1}}}\Big(\varphi\big(\levb[s_i][a](\lin,\btr[a]),\tau\big),\tb[s_i][a] + \tau\Big) \\
    \big(\levb[s_{i+1}][a,*](\lin), \tb[s_{i+1}][a,*]\big) &= \varsigma_{e_{s_i,s_{i+1}},\delta_{s_i,s_{i+1}}}\Big(\varphi\big(\levb[s_i][a,*](\lin),\overline\tau\big),\tb[s_i][a,*] + \overline\tau\Big).
\end{align*}
Lemma~\ref{lem:partialOrderAndVarphi}.4 and the iteration hypothesis give $\big(\levb[s_{i+1}][a](\lin,\btr[a]),\tb[s_{i+1}][a]\big)  \preceq \big(\levb[s_{i+1}][a,*](\lin),\tb[s_{i+1}][a,*]\big)$.
Hence, \eqref{eq:optimalSchedulingProofInductionHypothesis} is true for $i$ from $1$ to $k$.
We therefore have:
\begin{equation}\label{eq:dominanceInSk}
    \big(\levb[s_k][a](\lin,\btr[a]),\tb[s_k][a]\big)  \preceq \big( \levb[s_k][a,*](\lin),\tb[s_k][a,*]\big).
\end{equation}

\noindent Suppose that $\te[s_k][a,*](\lin)< \tb[s_k][a,*](\lin)$.
\begin{itemize}
    \item If $\lev[s_k][a,*](\lin) > 0$, by definition of $\btr[a,*](\lin)$, we would have $\te[v] + \sum_{i=1}^k \ti[s_{i-1},s_i] = \tb[s_k][a,*](\lin) <  \te[s_k][a,*](\lin) = \tb[v] - \fu[s_k,v]$, which means that~\eqref{eq:stationSequenceTimes} and therefore~\eqref{eq:optimalScheduilngProblem} do not admit solutions.
    \item If $\lev[s_k][a,*](\lin) = 0$, and since $\btr[a]$ is a solution of~\eqref{eq:stationSequenceTimes}, we have $\tb[s_k][a] \leq \te[s_k][a] = \tb[v] - \ti[s_k,v] = \te[s_k][a,*] < \tb[s_k][a,*]$. Hence~\eqref{eq:dominanceInSk} implies that $\levb[s_k][a](\lin,\btr[a]) = -\infty$, and therefore $\lev[v][a](\lin,\btr[a]) = -\infty$.
\end{itemize}

\noindent Suppose now that $\te[s_k][a,*](\lin)\geq \tb[s_k][a,*](\lin)$.
Since $\btr[a]$ satisfies~\eqref{eq:stationSequenceTimes}, and by definition of~$\btr[a,*](\lin)$, we have $\te[s_k][a] = \te[s_k][a,*](\lin) = \tb[v]- \ti[u,v]$.
Equation~\eqref{eq:optimalSchedulingProofInductionHypothesis} for $i=k$ and Lemma~\ref{lem:partialOrderAndVarphi}.5 then give $\big(\lev[s_{k}][a](\lin,\btr[a]),\te[s_{k}][a]\big)  \preceq \big(\lev[s_{k}][a,*](\lin),\te[s_{k}][a,*](\lin)\big)$, and therefore $\lev[s_{k}][a](\lin,\btr[a]) \leq \lev[s_{k}][a,*](\lin)$. As a consequence:
$$ \lev[v][a](\lin,\btr[a]) = \pos(\lev[s_{k}][a](\lin,\btr[a]) - \fu[s_k,v] - \fu[v]) \leq \pos(\lev[s_{k}][a,*](\lin) - \fu[s_k,v] - \fu[v]) = \lev[v][a,*](\lin).$$
This disjunction of cases concludes the proof of the theorem.
\halmos
\endproof

\subsection{Proof of Theorem~\ref{prop:stationSequenceDominance} on station sequences dominance.}
The first point of Theorem~\ref{prop:stationSequenceDominance} follows immediately from the recursive definition of $\btr[a,*]$.

Let us now consider the second point. 
Let $a$ and $a'$ be the two station sequences $u=s_0,s_1,\ldots,s_k, s_{k} = v$ and $u=s'_0,s'_1,\ldots,s'_{k'}, s'_{k'+1} = v$ and integers $i$ and $i'$ such that $k-i = k'-i'$ and $s_{i + j} = s'_{i' + j}$ for $j$ in $\{0,\ldots,k-i\}$. Suppose that:
$$\sum_{j=1}^i \co[s_{j-1},s_j] \leq \sum_{j=1}^{i'} \co[s'_{j-1},s'_j] 
\quad \text{and} \quad
\big(\levb[s_i][a,*](\lin),\tb[s_i][a,*](\lin)\big) 
\succeq 
\big(\levb[s'_{i'}][a',*](\lin),\tb[s'_{i'}][a,*](\lin)\big) \text{ for all } \lin \in [0,\mc].$$
Let us now prove by induction on $\tilde\imath$ from $i$ to $k$ that:
$$\sum_{j=1}^{\tilde\imath} \co[s_{j-1},s_j] \leq \sum_{j=1}^{i' + \tilde\imath - i} \co[s'_{j-1},s'_j] \text{ and }
\big(\levb[s_{\tilde\imath}][a,*](\lin),\tb[s_{\tilde\imath}][a,*](\lin)\big) 
\succeq 
\big(\levb[s'_{i' + \tilde\imath - i}][a',*](\lin),\tb[s'_{i' + \tilde\imath - i}][a,*](\lin)\big) \text{ for all } \lin \in [0,\mc].$$
The result for $\tilde\imath$ comes from the hypothesis.
Let us now suppose that the result is true up to $\tilde\imath < k$.
Since $s_{\tilde\imath} = s'_{i' + \tilde\imath - i}$ and $s_{\tilde\imath} = s'_{i' + \tilde\imath - i}$, we immediately obtain $\sum_{j=1}^{\tilde\imath + 1} \co[s_{j-1},s_j] \leq \sum_{j=1}^{i' + \tilde\imath - i + 1} \co[s'_{j-1},s'_j]$.
Besides, Lemma~\ref{lem:partialOrderAndVarphi}.4 applied with $\ell = \levb[s'_{i' + \tilde\imath - i}][a',*](\lin)$, $t = \tb[s'_{i' + \tilde\imath - i}][a,*](\lin)$, $\tau = \te[s'_{i' + \tilde\imath - i}][a,*](\lin) - \tb[s'_{i' + \tilde\imath - i}][a,*](\lin)$, 
$\overline l = \levb[s_{\tilde \imath}][a,*](\lin)$, $\overline t = \tb[s_{\tilde \imath}][a,*](\lin)$, $e = \fu[s_{\tilde\imath}, s_{\tilde\imath + 1}]$, and $\delta = \ti[s_{\tilde\imath}, s_{\tilde\imath + 1}]$, we get:
$$\big(\levb[s_{\tilde\imath + 1}][a,*](\lin),\tb[s_{\tilde\imath + 1}][a,*](\lin)\big) 
\succeq \big(\levb[s'_{i' + \tilde\imath - i + 1}][a',*](\lin),\tb[s'_{i' + \tilde\imath - i + 1}][a,*](\lin)\big) \text{ for all } \lin \in [0,\mc],$$
which give the induction and the result for $\tilde i = k$, i.e.,
$$\sum_{j=1}^k \co[s_{j-1},s_j] \leq \sum_{j=1}^{k'} \co[s'_{j-1},s'_j] 
\quad \text{and} \quad
\big(\levb[s_k][a,*](\lin),\tb[s_k][a,*](\lin)\big) 
\succeq 
\big(\levb[s'_{k'}][a',*](\lin),\tb[s'_{k'}][a,*](\lin)\big) \text{ for all } \lin \in [0,\mc].$$
Since $s_k = s'_{k'}$, we have $\te[s_k][a,*] = \te[s'_{k'}][a',*]$, and Lemma~\ref{lem:partialOrderAndVarphi}.5 gives $\big(\lev[s_k][a,*](\lin),\te[s_k][a,*](\lin)\big) 
\succeq 
\big(\lev[s'_{k'}][a',*](\lin),\te[s'_{k'}][a',*](\lin)\big) \text{ for all } \lin \in [0,\mc]$, and hence $\fc{a}(\lin) = \lev[s_k][a,*](\lin) - \fu[s_k,v] - \fu[v] \geq \lev[s'_{k'}][a',*](\lin) - \fu[s'_{k'},v] - \fu[v] = \fc{a'}(\lin)$, which gives the result.

\subsection{Proof of Theorem~\ref{thm:stationSequenceFu} and Corollary~\ref{cor:linearDominance} on Linear Recharge}

Let us start with the following corollary of Theorem~\ref{prop:optimalSchedulingOfStationSequence}.
\begin{corollary}\label{cor:a_feasibility_conditions}
    Given a station sequence $u=s_0,s_1,\ldots,s_k,s_{k+1}=v$, we have $\lev[v][a,*](\lin) = \fc{a}(\lin) > -\infty$ if and only if the three following conditions are satisfied:
    \begin{subequations}
        \label{eq:a_feasibility_conditions}
        \begin{align}
            \tb[s_k][a,*](\lin) &\leq \te[s_k][a,*](\lin) \label{eq:btr_feasibility} \\
            \lin & \geq e_{u,s_1} \label{eq:s1_reachable} \\
            \lev[s_k][a,*](\lin) &\geq e_{s_k,v} + e_v. \label{eq:v_reachable}
        \end{align}
    \end{subequations}
\end{corollary}
\proof{Proof.}
It follows from Theorem~\ref{prop:optimalSchedulingOfStationSequence} and the structure of $\btr[a,*](\lin)$. Indeed,
if~\eqref{eq:btr_feasibility} is not satisfied, then $\btr[a,*](\lin)$ is not feasible. 
Supposing that~\eqref{eq:btr_feasibility} is satisfied,
for any $i$ in $\{1,\ldots,k\}$, the level
$\lev[s_i][a](\lin, \btr[a,*])$ is non-smaller than $0$ if~\eqref{eq:s1_reachable} is satisfied and equal to $-\infty$ otherwise.
Finally, supposing that~\eqref{eq:btr_feasibility} and~\eqref{eq:s1_reachable} are satisfied, the level $\lev[v][a](\lin, \btr[a,*])$ is non-smaller than $0$ if~\eqref{eq:v_reachable} is satisfied and equal to $-\infty$ otherwise.
\halmos\endproof

\vspace*{1em}
\noindent
\proof{Proof of Theorem~\ref{thm:stationSequenceFu}.}
Let us place ourselves in the linear case ($\varphi$ given by Equation~\eqref{eq:linear_phi}).
We recall that $\fc{a}(\lin) > -\infty$ if and only if the three conditions of equation~\eqref{eq:a_feasibility_conditions} are satisfied.
Our proof is as follows: We start by proving that joint satisfaction of the three conditions of Equation~\eqref{eq:a_feasibility_conditions} is equivalent to a simple closed formula.
We then deduce from it the value of $\bc{a}(0)$, and then the other closed formula given in the proposition.

\vspace*{1em}
\noindent
\textbf{Simpler formula for $\tb[s_i][a,*](\lin)$ and $\lev[s_i][a](\lin, \btr[a,*])$.}
Given $e$ and $\ell$ in $[0,\mc]$ we have
$\tch(\ell,e) = \frac{e-\ell}{\alpha}.$
Using the definition of $\btr[a,*]$, an immediate induction on $i$ then gives that:
\begin{align}
&\tb[s_i][a,*](\lin) = \te[u] +  \underbrace{\sum_{j=1}^i \delta_{s_{j-1},s_j}}_{\substack{\text{Driving time required} \\ \text{to reach $s_i$}}} + \underbrace{\max \Big(0, \frac{1}{\alpha}\big(\overbrace{\sum_{j=1}^ie_{s_{j-1}s_j}}^{\substack{\text{Energy required}\\ \text{to reach $s_i$}}} - \lin \big) \Big)}_{\text{Charging time required to reach $s_i$}} && \text{for $i \in \{1,\ldots,k\}$} \\
&\lev[s_i][a](\lin, \btr[a,*]) = \max\Big(e_{s_j,s_{j+1}},\lin -\sum_{j=1}^i e_{s_{j-1}s_j}\Big) && \text{for $i \in \{1,\ldots,k-1\}$}. \label{eq:lev_si_lin}
\end{align}
\textbf{A simpler equivalent to~\eqref{eq:btr_feasibility}.}
Let us first prove that the following equation is equivalent to~\eqref{eq:btr_feasibility}.
\begin{equation}
    \label{eq:sufficientCOnditionBtrFeasibility}
    \lin \geq \sum_{j=1}^{k} e_{s_{j-1}s_j} - \underbrace{\alpha \Big(\tb[v] - \te[u] - \sum_{j=1}^{k+1}\delta_{s_{j-1},s_j}\Big)}_{\lav}
\end{equation}
Remark that \eqref{eq:btr_feasibility} is therefore equivalent to:
\begin{equation}
    \label{eq:firstEquaivalentToBtrFeas}
    \tb[v] - \ti[s_k,v] \geq \te[u] +  \sum_{j=1}^k \delta_{s_{j-1},s_j} + \max \Big(0, \frac{1}{\alpha}\sum_{j=1}^k\big( e_{s_{j-1}s_j} - \lin \big) \Big).
\end{equation}
If the maximum on the right-hand side of~\eqref{eq:firstEquaivalentToBtrFeas} is $0$,~\eqref{eq:stationSequenceTime} ensures that both \eqref{eq:firstEquaivalentToBtrFeas} and \eqref{eq:sufficientCOnditionBtrFeasibility} are satisfied. 
Otherwise, it is immediate that \eqref{eq:firstEquaivalentToBtrFeas} is satisfied if and only if \eqref{eq:sufficientCOnditionBtrFeasibility} is satisfied.

\vspace*{1em}
\noindent
\textbf{A formula for $\lev[s_k][a](\lin, \btr[a,*])$ when \eqref{eq:btr_feasibility} is satisfied.}
Supposing that~\eqref{eq:sufficientCOnditionBtrFeasibility} is satisfied, using the definition of $\varphi$, the fact that $\te[s_k][a,*](\lin)  = \tb[v] - \ti[s_k,v]$, and the result above, we have:
\begin{equation}\label{eq:lev_a_sk}
    \begin{array}{rl}
& \lev[s_k][a](\lin, \btr[a,*]) = \max \bigg[\mc, \\
& \quad 
\underbrace{\max\Big(0,\lin - \sum_{j=1}^k e_{s_{j-1}s_j}\Big)}_{\lev[s_k][a](\lin, \btr[a,*]) - e_{s_{k-1},s_k}} 
+ \alpha
\Big(\underbrace{\tb[v] - \ti[s_k,v]}_{\te[s_k][a,*]} - \underbrace{\Big[\te[u] +  \sum_{j=1}^k \delta_{s_{j-1},s_j} + \max \Big(0, \frac{1}{\alpha}\big(e_u+\sum_{j=1}^k e_{s_{j-1}s_j} - \lin \big) \Big)  \Big]}_{\te[s_k][a,*]} \Big) \bigg] \\
&=
\max\bigg[\mc,
\lin -\sum_{j=1}^k e_{s_{j-1}s_j} + \alpha\Big(\tb[v] - \ti[s_k,v] -  \Big[\te[u] +  \sum_{j=1}^k \delta_{s_{j-1},s_j}\Big] \Big) \bigg].
    \end{array}
\end{equation}

\noindent
\textbf{A condition equivalent to~\eqref{eq:v_reachable} and \eqref{eq:btr_feasibility}.}
Given~\eqref{eq:stationSequenceTime}, $\mc \geq e_{s_k,v}+e_v$. We therefore deduce from~\eqref{eq:lev_a_sk} that, provided that~\eqref{eq:btr_feasibility} is satisfied, Condition~\eqref{eq:v_reachable} is equivalent to: 
\begin{equation} \label{eq:v_reachableAsAFunctionOfLin}
     \lin \geq \sum_{j=1}^{k+1} e_{s_{j-1}s_j} +e_v - \underbrace{\alpha \Big(\tb[v] - \te[u] - \sum_{j=1}^{k+1}\delta_{s_{j-1},s_j}\Big)}_{\lav} = - \ldinc.
\end{equation}
Remark that~\eqref{eq:v_reachableAsAFunctionOfLin} implies the sufficient condition~\eqref{eq:sufficientCOnditionBtrFeasibility} for~\eqref{eq:btr_feasibility}. 
We have proved that:
\begin{itemize}
    \item If~\eqref{eq:btr_feasibility} is satisfied, Condition~\eqref{eq:v_reachable} is equivalent to~\eqref{eq:v_reachableAsAFunctionOfLin},
    \item and~\eqref{eq:v_reachableAsAFunctionOfLin} implies \eqref{eq:btr_feasibility}.
\end{itemize}
Hence~\eqref{eq:v_reachableAsAFunctionOfLin} is equivalent to $\big(\text{\eqref{eq:v_reachable} and \eqref{eq:btr_feasibility}}\big)$.

\vspace*{1em}
\noindent
\textbf{Proof of the closed formula for $\lmb$.}
Corollary~\eqref{cor:a_feasibility_conditions} and the previous equivalence therefore imply that $\fc{a}(\lin) > -\infty$ if and only if~\eqref{eq:s1_reachable} and~\eqref{eq:v_reachableAsAFunctionOfLin}, that is, if and only if:
$$\lin \geq \max(e_{u} + e_{u,s_1},- \ldinc).$$ 
We therefore have:
$$\lmb = \max(e_{u} + e_{u,s_1},- \ldinc).$$

\noindent
\textbf{Proof of the closed formula for $\fc{a}(\lin)$ and $\lme$.}
Supposing $\lin \geq \lmb$, and using~\eqref{eq:lev_a_sk}, we have:
\begin{equation}
    \label{eq:fac_lin_detailed}
\begin{array}{rl}
\fc{a}(\lin) = \lev[v][a](\lin, \btr[a,*]) &= \lev[s_k][a](\lin, \btr[a,*]) - e_{s_k,v} - e_v  \\
&= \min\bigg(\mc, \lin + \underbrace{\alpha \Big(\tb[v] - \te[u] - \sum_{j=1}^{k+1}\delta_{s_{j-1},s_j}\Big)}_{\lav} - \displaystyle\sum_{i=1}^{k}\fu[s_{i-1},s_i]\bigg) - e_{s_k,v} -e_v
\end{array}
\end{equation}
This result with $\lin = \mc$ immediately gives:
\begin{align*}
    \lme &=
\min\bigg(\mc, \mc+\lav - \displaystyle\sum_{i=1}^{k}\fu[s_{i-1},s_i]\bigg) - \fu[s_k,v] - e_v.
\end{align*}

\noindent
\textbf{Proof of the closed formula for $\ldi$.}
Considering~\eqref{eq:fac_lin_detailed} with $\lin = \bc{a}(0)$, we have two cases: either $\lev[s_k][a](\lin, \btr[a,*]) = \mc$, and we obtain $\fc{a}(\lmb) = \fc{a}(\lme)$. Or $\lev[s_k][a](\lin, \btr[a,*]) \leq \mc$, and $\fc{a}(\lmb) = \lmb + \ldinc$. We therefore get:
$$\ldi = \min\big(\ldinc, \lme -\lmb \big).$$

\noindent
\textbf{Proof of the closed formula for $\fc{a}(\lin)$ and $\bc{a}(\lout)$.}
Considering again~\eqref{eq:fac_lin_detailed} for a general $\lin$ and the two cases $\lev[s_k][a](\lin, \btr[a,*]) = \mc$ and $\lev[s_k][a](\lin, \btr[a,*]) \leq \mc$, and substituting the values of $\lme$ and $\ldi$, we obtain:
\begin{alignat*}{1}
\fc{a}(\lin) &= 
\begin{cases}
\min\big(\lme,\lin + \ldi\big) & \text{if } \lin \geq \lmb, \\
-\infty & \text{otherwise.}
\end{cases}
\end{alignat*}
Getting back to the definition of $\bc{a}$ in Equation~\eqref{eq:definition_of_bca}, we get:
\begin{alignat*}{1}
\bc{a}(\lout) &= 
\begin{cases}
\lmb + \max(0, \lout - \ldi) & \text{if } \lout \leq \lme, \\
+\infty & \text{otherwise.} 
\end{cases}
\end{alignat*}
This concludes the proof of Theorem~\ref{thm:stationSequenceFu}.
\halmos\endproof

\proof{Proof of Corollary~\ref{cor:linearDominance}.}
Let $u=s_0,\ldots,s_k=v$ be a station sequence, and $\lin$ be in $[0,M]$. Then Theorem~\ref{thm:stationSequenceFu} ensures that:
\begin{equation}
\label{eq:linearBeginChargeAndTime}
\begin{array}{ll}
    \begin{cases}
        \levb[s_i][a,*] = \max(0,\lin - \sum_{j=1}^i \fu[s_{j-1},s_j])  \\
        \tb[s_i][a,*] = \te[v] + \sum_{j=1}^{i} \ti[s_{i-1},s_i] + \frac{1}{\alpha}\max\Big(0,\sum_{j=1}^i \fu[s_{i-1},s_i] - \lin\Big)
    \end{cases}
    & \text{if } \lin \geq \fu[u,s_1] \\
    \begin{cases}
        \levb[s_i][a,*] = -\infty  \\
        \tb[s_i][a,*] = +\infty
    \end{cases}
    & \text{otherwise.} 
\end{array}
\end{equation}
Let us now place ourselves in the setting of Corollary~\ref{cor:linearDominance}.
If $\fu[u,s_1] \leq \fu[u,s_2]$, Equation~\eqref{eq:linearBeginChargeAndTime} gives
$\big(\levb[s_i][a,*],\tb[s_i][a,*](\lin)\big) 
\succeq 
\big(\levb[s'_{i'}][a',*],\tb[s'_{i'}][a,*](\lin)\big)$ for $\lin \leq \fu[u,s_2]$.
Furthermore, $\big(\levb[s_i][a,*],\tb[s_i][a,*](\lin)\big) 
\succ 
\big(\levb[s'_{i'}][a',*],\tb[s'_{i'}][a,*](\lin)\big)$ 
for $\fu_{u,s_1}<\lin \leq \fu[u,s_2]$.

Consider now $\lin \geq \max(\fu_{u,s_1}, \fu[u,s_2])$.
Recall that $\tch(\ell,\ell') = \frac{\ell'-\ell}{\alpha}.$
Hence, using Equation~\eqref{eq:linearBeginChargeAndTime}, we get:
\begin{align*}
    \tb[s_i][a,*] + \tch\big( \levb[s_i][a,*],  \levb[s'_{i'}][a',*]\big) 
    &= \te[v] + \sum_{j=1}^{i} \ti[s_{j-1},s_j] + \frac{1}{\alpha}\max\Big(0,\sum_{j=1}^i \fu[s_{j-1},s_j] - \lin\Big) \\
    & \quad + \frac{1}{\alpha}\Big(\max(0,\lin - \sum_{j=1}^{i'} \fu[s'_{j-1},s_j]) - \max(0,\lin - \sum_{j=1}^i \fu[s_{j-1},s_j])  \Big) \\
    &= \te[v] + \sum_{j=1}^{i} \ti[s_{j-1},s_j] + \frac{1}{\alpha}\sum_{j=1}^i \fu[s_{j-1},s_j] - \lin + \frac{1}{\alpha}\Big(\max(0,\lin - \sum_{j=1}^{i'} \fu[s'_{j-1},s_j]).
\end{align*}
We have $\big(\levb[s_i][a,*],\tb[s_i][a,*](\lin)\big) 
\succeq 
\big(\levb[s'_{i'}][a',*],\tb[s'_{i'}][a,*](\lin)\big)$  if and only if $\tb[s_i][a,*] + \tch\big( \levb[s_i][a,*],  \levb[s'_{i'}][a',*]\big) \leq \tb[s'_{i'}][a',*]$, that is,
if and only if:
\begin{align*}
   \te[v] + \sum_{j=1}^{i} \ti[s_{j-1},s_j] + \frac{1}{\alpha}\sum_{j=1}^i \fu[s_{j-1},s_j] - \lin + \frac{1}{\alpha}\max\Big(0,\lin - \sum_{j=1}^{i'} \fu[s'_{j-1},s_j]\Big) \\
   \quad  \leq \te[v] + \sum_{j=1}^{i'} \ti[s'_{j-1},s_j] + \frac{1}{\alpha}\max\Big(0,\sum_{j=1}^{i'} \fu[s_{j-1},s_j] - \lin\Big),
\end{align*}
that is, if and only if:
$$ \sum_{j=1}^{i} \ti[s_{j-1},s_j] + \frac{1}{\alpha}\sum_{j=1}^i \fu[s_{j-1},s_j] \leq  \sum_{j=1}^{i'} \ti[s'_{j-1},s_j] + \frac{1}{\alpha}\sum_{j=1}^{i'} \fu[s_{j-1},s_j],$$
which concludes the proof of Corollary~\ref{cor:linearDominance}.
\halmos\endproof

\end{document}